\documentclass[11pt,english,a4paper,reqno]{amsart}
\usepackage[%
   pdftitle={Variational characterizations of the effective multiplicative factor},
   pdfauthor={B. Lods},
   pdfsubject={Variational characterizations of the effective multiplicative factor},
  colorlinks=true,
  linkcolor=blue,
  citecolor=blue,
]{hyperref}
\usepackage[T1]{fontenc}
\usepackage{mathrsfs}
\usepackage{amsmath,amsthm,amssymb}
\usepackage{latexsym}
\usepackage{times}
\usepackage{enumerate}
\usepackage{mathrsfs}
\usepackage{stmaryrd}
\usepackage{amsopn}
\usepackage{amsmath}
\usepackage{amssymb}
\usepackage{amsfonts}
\usepackage{amsbsy}
\usepackage{amscd,indentfirst}
\usepackage{hyperref}
\newtheorem{theorem}{\rm \bf Theorem}[section]

\newtheorem{proposition}[theorem]{Proposition}
\newtheorem{corollary}[theorem]{Corollary}

\newtheorem{definition}[theorem]{Definition}

\newtheorem{remark}[theorem]{Remark}
\font\teneuf=eufm10 at 12pt \font\seveneuf=eufm7 at 8pt
\font\fiveeuf=eufm5 at 6pt
\newfam\euffam
\textfont\euffam=\teneuf \scriptfont\euffam=\seveneuf
\scriptscriptfont\euffam=\fiveeuf

\def\bea{\begin{eqnarray}}
\def\eea{\end{eqnarray}}
\def\beao{\begin{eqnarray*}}
\def\eeao{\end{eqnarray*}}

\numberwithin{equation}{section}
\newfont{\secgoth}{eufm10 at 16pt}

\def\di {\text{div}_x}
\def \esup {\operatornamewithlimits{ess\,sup}}
\def \esinf {\operatornamewithlimits{ess\,inf}}
\def \ke {k_{\mathrm{eff}}}
\def \d {{d}}
\def \ge {k_{\mathrm{eff}}}
\def \geq {\geqslant}
\def \leq {\leqslant}
\def \kats {\kappa_s}
\def \katf {\kappa_f}
\def \katsdo {\kappa_s(\cdot,\cdot,\cdot)}
\def \katfdo {\kappa_f(\cdot,\cdot,\cdot)}
\def \ds {\displaystyle}

\title[Variational characterizations of the effective multiplication factor]
{Variational characterizations of the effective multiplication
factor of a nuclear reactor core}

\author[B. \textsc{ Lods}]{}
\begin{document}
\bibliographystyle{plain}
\maketitle \centerline{{B}ertrand {\sc Lods}} \smallskip {
\centerline{Universit\'{e} Blaise Pascal (Clermont II)}
   \centerline{Laboratoire de Math\'ematiques, CNRS UMR 6620}
   \centerline{63117 Aubi\`{e}re, France.}
   \centerline{\texttt{lods@math.univ-bpclermont.fr}}
   }

\begin{abstract} We prove some inf--sup and sup--inf formulae for the so--called
effective multiplication factor arising in the study of reactor
analysis. We treat in a same formalism the transport equation and
the energy--dependent diffusion equation.


\end{abstract}
%
%
%
\section{Introduction}

The aim of this paper is to give some variational characterizations
of the {\it effective multiplication factor} arising in nuclear
reactor theory. This work follows a very recent paper by M.
Mokhtar-Kharroubi \cite{mk} devoted to the leading eigenvalue of
transport operators.

In practical situations, the \textit{power distribution} in a {
stable} nuclear reactor core   is determined as the steady-state
solution $\phi$ of a linear transport equation for the neutron flux.
Because of interactions between neutrons and fissile isotopes, a
fission chain reaction occurs in the reactor core. Precisely, when
an atom undergoes nuclear fission, some neutrons  are ejected from
the reaction and subsequently shall interact with the surrounding
medium. If more fissile fuel is present, some may be absorbed and
cause more fissions (see for details
Refs.~\cite{bell,birk1,dud,planchard}). The linear stationary
transport equation is therefore of non-standard type in the sense
that the source term is itself a function of the solution. When
delayed neutrons are neglected, this equation reads
\begin{multline}
\label{equ} v\cdot \nabla_x \phi(x,v)+
\sigma(x,v)\phi(x,v)-\int_{V}\kats(x,v,v')\phi(x,v')\d\mu(v')=\\
\frac{1}{\ge}\int_{V}\katf(x,v,v')\phi(x,v')\d\mu(v')
\end{multline}
with {\it free-surface} boundary condition (i.e. the incoming flux
is null). Here, the unknown $\phi(x,v)$ is the neutron density at
point $x \in \mathcal{D}$ and velocity $v \in V,$ where
$\mathcal{D}$ is an open subset of $\mathbb{R}^N$ (representing the
reactor core) and the velocity space $V$ is a closed subset of
$\mathbb{R}^N$, $\d\mu(\cdot)$ being a positive Radon measure
supported by $V$. For the usual cases, $\d\mu(\cdot)$ is either the
Lebesgue measure on $\mathbb{R}^N$ (continuous model) or on spheres
(multigroup model). The  transfer cross-sections $\katsdo$ and
$\katfdo$ describe respectively the {\it pure scattering} and the
{\it fission} process. The {\it nonnegative} bounded function
$\sigma(\cdot,\cdot)$ is the {\it absorption cross--section}
\cite{bell,dud,planchard}. The positive ratio $\ge$ is called the
{\it criticality eigenvalue} (or the {\it effective multiplication
factor}). It represents the average  number of neutrons that go on
to cause another fission reaction. The remaining neutrons either
fail to induce fission, or  never get absorbed and exit the system.
Consequently, $\ge$ measures the balance between the number of
neutrons in successive generations (where the birth event separating
generations is the fission process). The interpretation of the
effective multiplication factor $\ge$ is related to the following
three cases
(see \cite{bell,dud,planchard}):\\

 \indent $\bullet$ If
$\ge=1$, there is a perfect balance between production and removal
of neutrons. The reactor is then said to be {\it critical}.

\indent $\bullet$ The reactor is \textit{sub-critical} when $\ge <
1$. This means that the removal of neutrons (at the boundary or due
to absorption by the surrounding media) excesses the fission process
and the chain reaction dies out rapidly.

\indent $\bullet$ When $\ge > 1$, the fission chain reaction grows
without
bound and the reactor is said to  be {\it super--critical}.\\

 Up to now, in practical applications,
the effective multiplication factor $\ge$ was usually given by
$r_{\sigma}\left[(\mathcal{T}-\mathcal{K}_{s})^{-1}
\mathcal{K}_{f}\right]=\ge,$   where the precise definition of the
operators $\mathcal{T},$ $\mathcal{K}_s$ and $\mathcal{K}_f$ is
given subsequently and $r_\sigma[B]$ denotes the spectral radius of
any generic bounded operator $B$: $r_\sigma[B]=\lim_{n \to
\infty}\|B^n\|^{1/n}$. Because its requires the computation of the
resolvent $(\mathcal{T}-\mathcal{K}_s)^{-1}$, such a
characterization makes the analysis of the effective multiplication
factor quite difficult to handle. In particular, practical estimates
of $\ge$ as well as computational approximations are rather involved
and merely rely on (direct or inverse) power method \cite{berry}.
Our main objective in this paper is to provide tractable variational
characterizations of $\ge$ in terms of the different data of the
system  and suitable test functions. We hope that such a
characterization shall be of interest for practical computations or
for the homogenization of the above criticality transport equation
in periodic media \cite{allaire,bal,hom}. We also believe that our
characterization can be useful in the delicate optimization problem
of the assembly distribution in a nuclear reactor (see the recent
contribution \cite{thevenot} based upon the homogenization method
and where the criticality eigenvalue $\ge$ plays a crucial role).

To be precise, we provide here {\it variational characterizations}
of the effective multiplication factor of the type

\begin{equation} \label{car}
\begin{split}
\frac{1}{\ge}&=\underset{\varphi \, \in W_p^+}{\min}\;
\underset{(x,v) \in \mathcal{D} \times V}{\esup} \frac{v\cdot
\nabla_x \varphi(x,v)+ \sigma(x,v)\varphi(x,v)-\ds
\int_{V}\Sigma_s(x,v,v')\varphi(x,v')\d\mu(v')}
{\ds \int_{V}\Sigma_f(x,v,v')\varphi(x,v')\d\mu(v')}\\
& \\
&=\underset{\varphi \, \in W_p^+}{\max}\; \underset{(x,v) \in
\mathcal{D} \times V}{\esinf} \frac{v\cdot \nabla_x \varphi(x,v)+
\sigma(x,v)\varphi(x,v)-\ds
\int_{V}\Sigma_s(x,v,v')\varphi(x,v')\d\mu(v')} {\ds
\int_{V}\Sigma_f(x,v,v')\varphi(x,v')\d\mu(v')},
\end{split}
\end{equation}
where $W_p^+$ is a suitable class of {\it positive test--functions}
in $L^p(\mathcal{D} \times V,\d x\d\mu(v))$ $(1 \leq p < \infty).$
This result (Theorem \ref{princitran}) holds true under compactness
assumption on the {\it full collision operator}
$$\mathcal{K}\::\:\psi 
\longmapsto \int_V
\left(\Sigma_s(x,v,v')+\Sigma_f(x,v,v')\right)\psi(x,v')\,\d\mu(v')
$$ and under {\it positivity
assumptions on the fission cross--section} $\katfdo$. The main
strategy to derive (\ref{car}) is adapted from \cite{mk} where M.
Mokhtar-Kharroubi proved similar variational characterizations for
the leading eigenvalue of {\it perturbed} transport operators. Note
that the above characterization still holds true for transport
equations with general  boundary conditions modeled by some
nonnegative albedo operator (see Remark \ref{boundary}).

At this point, one recalls that, besides the critical eigenvalue
$\ge,$  it is also possible to investigate the reactivity of the
nuclear reactor core through another physical parameter, namely, the
leading eigenvalue $s(\mathcal{A})$ of the operator
$\mathcal{A}=\mathcal{T}+\mathcal{K}_s+\mathcal{K}_f $, also
associated to positive eigenfunctions. The two parameters $\ge$ and
$s(\mathcal{A})$ are related by the following: if $s(\mathcal{A}) <
0$ then the reactor is subcritical (i.e. $\ge <1$), while it is
super-critical whenever $s(\mathcal{A})
>0$. The reactor is critical when $s(\mathcal{A})=0$. The paper \cite{mk} provides a variational characterization of the leading
eigenvalue of the transport operator $\mathcal{A}$. However, for
practical calculations in nuclear engineering, the critical
eigenvalue $\ge$ is a more effective parameter. Actually, the
existence of the leading eigenvalue
$s(\mathcal{T}+\mathcal{K}_s+\mathcal{K}_f)$ is not always ensured
but is related to the size of the domain $\mathcal{D}$ and the
possibility of small velocities (for more details on this
disappearance phenomenon, see e.g. \cite[Chapter 5]{mmk}). Since the
existence of the effective multiplication factor is not restricted
by such physical constraints, it appears more efficient to measure
the reactivity of nuclear reactor cores by $\ge$. This is what
motivated us to generalize the result of \cite{mk} and provide
variational characterization of $\ge$.
\smallskip

In this paper, we also give a characterization of the criticality
eigenvalue associated to the \textit{energy-dependent diffusion}
model  used in nuclear reactor theory \cite{bell,pao2,planchard}.
For this description, the critical problem reads
\begin{multline}\label{diffu1}
-\di(\mathbb{D}(x,\xi)\nabla_x \varrho(x,\xi))+\sigma(x,\xi)\varrho(x,\xi)-\int_{E}\Sigma_s(x,\xi,\xi')\,\varrho(x,\xi')\,\d\xi'\\
=\frac{1}{\ke} \int_E \Sigma_f(x,\xi,\xi')\,\varrho(x,\xi')\,\d
\xi',\end{multline} complemented by the Dirichlet boundary
conditions $\varrho(\cdot,\xi)_{|\partial \mathcal{D}}=0$ a.e. $\xi
\in E.$ Here $E$ is the set of admissible energies $\xi=\frac{1}{2}m
v^2$ ($m$ being the neutron mass and $v$ the velocity), i.e. $E$ is
a subset of $[0,+\infty[$. The diffusion coefficient
$\mathbb{D}(\cdot,\cdot)$ is a matrix--valued function over
$\mathcal{D} \times E$ and the unknown $\varrho(\cdot,\cdot)$ is
\textit{nonnegative}.

The derivation of diffusion-like models for some macroscopic
distribution function $\varrho(x,\xi)$ (corresponding to some
angular momentum of the solution $\phi$ to \eqref{equ}) is motivated
in nuclear engineering by the necessity to provide simplified models
tractable numerically. Such a \textit{energy-dependent diffusion}
model can be derived directly from a phenomenological analysis of
the scattering models or it can be derived from the above kinetic
equation \eqref{equ} through a suitable asymptotic procedure (see
\cite{degond} and the references therein for more details on that
matter).

For this \textit{energy-dependent diffusion} model, we give a
variational characterization of $\ke$ in terms of sup-inf and
inf-sup criteria in the spirit of \eqref{car}.\medskip

Actually, to treat the two above problems \eqref{equ} and
\eqref{diffu1} it is possible to adopt a unified mathematical
formalism. Precisely, let us denote by $\mathcal{K}_{s}$ the
integral operator with kernel $\katsdo$ and denote by
$\mathcal{K}_{f}$ the integral operator with kernel $\katfdo$. Then,
problems \eqref{equ} and \eqref{diffu1} may be written in a
\textit{unified abstract way}:
$$\mathcal{T}\phi_\mathrm{eff}+\mathcal{K}_{s}\phi_\mathrm{eff}+\frac{1}{\ge}\mathcal{K}_{f}\phi_\mathrm{eff}=0,\qquad \qquad \phi_\mathrm{eff} \geq 0,$$
where the unbounded operator $\,\mathcal{T}$ refers to, according to
the model we adopt:

\indent $\bullet$ the transport operator:
$$\mathcal{T}\phi(x,v)=-v \cdot \nabla_x \phi(x,v) -\sigma(x,v)\phi(x,v),$$
associated to the absorbing boundary conditions
$\phi_{|\Gamma_-}=0.$

\indent $\bullet$ the energy--dependent diffusion operator:
$$\mathcal{T}\psi(x,\xi)=\di\left(\mathbb{D}(x,\xi)\nabla_x\psi(x,\xi)\right)-\sigma(x,\xi)\psi(x,\xi)$$
with Dirichlet boundary conditions $\psi_{|\,\partial
\mathcal{D}}(\cdot,\xi)=0$ $(\text{a.e. } \xi \in E).$\\

The abstract treatment of the above problem is performed in Section
\ref{caraabs} and relies mainly on positivity and compactness
arguments. The main \textit{abstract} result of this paper (Theorem
\ref{supinfabs}) characterizes the criticality eigenvalue of a large
class of (abstract) unbounded operators in $L^p$-spaces. Besides
this main analytical result, we also prove abstract results with
their own interest. In particular, we provide in Theorem
\ref{theo:app} \textit{an approximation resolution for the
criticality eigenfunction} $\phi_\mathrm{eff}$ which shall be
hopefully useful for practical numerical approximations.

The outline of the paper is as follows. In Section 2, we describe
the unified and abstract framework which allows us to treat in a
same formalism Problems \eqref{equ} and \eqref{diffu1} with the aim
of establishing general inf--sup and sup--inf formulae for the
criticality eigenvalue of a class of unbounded operator.  In Section
3, we are concerned with the characterization of the effective
multiplication factor $\ge$ associated to the transport problem
\eqref{equ}. In Section 4, we investigate the effective
multiplication factor associated to the energy-dependent diffusion
model \eqref{diffu1}.
%
%
%
%
%
%
%

\section{Abstract variational characterization}\label{caraabs}
%
%
%
%
This section is devoted to several abstract variational
characterizations of the criticality eigenvalue. It is this abstract
material that shall allow us to treat in the same formalism Problems
\eqref{equ} and \eqref{diffu1}.

\subsection{Setting of the problem and existence result} Let us introduce the functional framework we shall use in the
sequel. Given a measure space $(\mathbf{\Omega},\nu)$ and a fixed $1
\leq p < \infty,$ define
$$X_p=L^p(\mathbf{\Omega},\d\nu)$$
and denote by $X_q$ its dual space, i.e.
$X_q=L^q(\mathbf{\Omega},\d\nu)$ $(1/p+1/q=1).$ We first recall
several definitions and facts about {\it positive operators}. Though
the various concepts we shall deal with could be defined in general
complex Banach lattices, we restrict ourselves to operators in $X_p$
$(1 \leq p < \infty)$:
\begin{definition}
A bounded operator $B$ in $X_p$ is said to be
\textit{\textbf{irreducible}} if, for every {\it nonnegative}
$\varphi \in X_p\setminus\{0\}$ and {\it nonnegative} $\psi  \in X_q
\setminus \{0\}$, there exists $n \in \mathbb{N}$ such that
$$\left \langle B^n\varphi,\psi \right \rangle > 0,$$
where $\left \langle \cdot,\cdot \right \rangle$ is the duality
pairing between $X_p$ and the dual space $X_q.$\end{definition} Let
us denote the set of {\it quasi-interiors} elements of $X_p$ by
$X_p^+,$ i.e.
$$X_p^+=\{f \in X_p\;;\;f(\omega) > 0\;\d\nu-\text{a.e. } \omega \in
\mathbf{\Omega}\}.$$
Notice that, if $f \in X_p^+$, then $\langle f, \psi \rangle
>0$ for any nonnegative $\psi \in X_q \setminus \{0\}$.
\begin{definition}
 A bounded operator $B$ in $X_p$ will be said to be
\textbf{\textit{positivity improving}} if its maps nonnegative $f
\in X_p \setminus \{0\}$ into $X_p^+$, i.e.
$$f \in X_p \setminus \{0\}, f \geq 0 \Longrightarrow B f \in
X_p^+.$$
\end{definition}
\begin{remark} Notice that, given a bounded operator $B$ in $X_p$,
if some power of $B$ is positivity improving, then $B$ is
irreducible. This provides a practical criterion of irreducibility.
\end{remark}
Recall also several fact about  power-compact  operators.
\begin{definition}
A bounded operator $B$ in a Banach space $X$ is said to be
\textit{\textbf{power-compact}} if there exists $n \in \mathbb{N}$
such that $B^n$ is a compact operator in $X$.
\end{definition}
The following fundamental result is due to B. De Pagter \cite{bdp}.
\begin{theorem}\label{bendp}
Let $B$ be a bounded operator in a Banach space $X$. If $B$ is
irreducible and power-compact then $r_{\sigma}(B) >0$ where
$r_{\sigma}(B)$ denotes the spectral radius of $B$.
\end{theorem}

 Let $\mathcal{T}$ be
a given  densely defined unbounded operator
$$\mathcal{T} \::\:\mathscr{D}(\mathcal{T}) \subset X_p \longrightarrow X_p$$
such that
\begin{equation}\label{assumptionT}
s(\mathcal{T}) < 0 \qquad \text{ and }\qquad
(0-\mathcal{T})^{-1}(X_p^+) \subset X_p^+.\end{equation} Let
$\mathcal{K}_{s}$ and $\mathcal{K}_{f}$ be two \textit{nonnegative
bounded} operators in $X_p$. We are interested in the abstract
critical problem:
\begin{equation}\label{abs}
(\mathcal{T}+\mathcal{K}_{s}+\frac{1}{\ge}\mathcal{K}_{f})\phi_\mathrm{eff}=0,
\qquad \phi_\mathrm{eff} \in
\mathscr{D}(\mathcal{T})\,,\,\phi_\mathrm{eff} \geq 0,
\:\phi_{\mathrm{eff}} \neq 0.
\end{equation} Since $s(\mathcal{T})<
0$, Problem \eqref{abs} is equivalent to
$$(0-\mathcal{T})^{-1}(\mathcal{K}_{s}+\frac{1}{\ge}\mathcal{K}_{f})\phi_\mathrm{eff}=\phi_\mathrm{eff}\,,\qquad
\phi_\mathrm{eff} \geq 0,\:\phi_{\mathrm{eff}} \neq 0.$$ Let us
introduce the family of operators indexed by the positive parameter
$\gamma$:
$$\mathcal{K}(\gamma)=\mathcal{K}_{s}+\frac{1}{\gamma}\mathcal{K}_{f}\qquad \qquad \gamma> 0.$$
Therefore, solving \eqref{abs} is equivalent to prove the existence
(and uniqueness) of $\ge > 0$ such that $1$ is an eigenvalue of
$(0-\mathcal{T})^{-1}\mathcal{K}(\ge)$ associated to a
\textit{nonnegative eigenfunction}. Such an existence and uniqueness
result can be found in \cite[Theorem 5.30]{mmk}   (see also
\cite{marek1}). We set $\mathcal{K}=\mathcal{K}_s+\mathcal{K}_f$
%
%
%
\begin{theorem}
\label{exist} Assume that $(0-\mathcal{T})^{-1}\mathcal{K}$ is
power-compact and that $(0-T)^{-1}\mathcal{K}_f$ is  irreducible.
Then, the spectral problem \eqref{abs} admits a unique solution $\ke
> 0$ associated with a nonnegative eigenfunction
$\phi_{\mathrm{eff}}$ if and only if
\begin{equation}
\label{iff} \lim_{\gamma \to
0}r_{\sigma}[(0-\mathcal{T})^{-1}\mathcal{K}(\gamma)] >1 \quad
\mbox{ and } \quad r_{\sigma}[(0-\mathcal{T})^{-1}\mathcal{K}_{s}] <
1.
\end{equation}
\end{theorem}

\begin{remark}\label{unifN} Notice that our assumptions differs slightly from that of
\cite{mmk}. Actually, in \cite{mmk}, it is assumed that
$(0-\mathcal{T})^{-1}\mathcal{K}(\gamma)$ is power-compact and
irreducible for any $\gamma >0$.  Our assumption implies those of
\cite{mmk}. Indeed, in this case, there is an integer $N \in
\mathbb{N}$ such that
$\left[(0-\mathcal{T})^{-1}\mathcal{K}\right]^N$ is compact.  Since,
for any $\gamma
>0$, $\mathcal{K}(\gamma) \leq \max\{1,1/\gamma\}\mathcal{K}$, one gets by a domination
argument that
$\left[(0-\mathcal{T})^{-1}\mathcal{K}(\gamma)\right]^N$ is compact
for any $\gamma
>0$. This means that the power at which $(0-\mathcal{T})^{-1}\mathcal{K}(\gamma)$
becomes compact is independent of $\gamma > 0$. In the same way,
since $(0-\mathcal{T})^{-1}\mathcal{K}(\gamma) \geq
\frac{1}{\gamma}(0-\mathcal{T})^{-1}\mathcal{K}_f$ for any $\gamma
>0$, our assumption implies the irreducibility of
$(0-\mathcal{T})^{-1}\mathcal{K}(\gamma)$ for any $\gamma >0.$
Notice also that the result still holds if
$(0-\mathcal{T})^{-1}\mathcal{K}_s$ is irreducible.
\end{remark}

\begin{remark}
\label{rem1} Under the assumptions of the previous Theorem, we point
out that the mapping $\gamma > 0 \mapsto
r_{\sigma}[(0-\mathcal{T})^{-1}\mathcal{K}(\gamma)]$ is continuous
(see \cite[Remark 3.3, p. 208]{kato}). By analyticity arguments
(Gohberg-Shmulyan theorem), it is also {\it strictly decreasing}.
Thus, $\ke$ is \textit{\textbf{characterized}} by
$$r_{\sigma}[(0-\mathcal{T})^{-1}\mathcal{K}(\ke)]=1.$$
\end{remark}

Let us now give some variational characterizations of the
criticality eigenvalue $\ke$ appearing in Theorem \ref{exist}.

\subsection{Abstract variational characterization of $\ke$} From now on, Assumption \eqref{iff} is assumed to be
fulfilled.
%
Let \begin{equation} \label{wp+} W_p^{+}:=\mathscr{D}(\mathcal{T})
\cap X^{+}_p.
\end{equation}
We start with the following characterization of $\ke$ in terms of
{\it super-solution} to the spectral problem \eqref{abs}.
\begin{proposition}
\label{tau} Assume that $(0-\mathcal{T})^{-1}\mathcal{K}(\gamma)$ is
power-compact and irreducible for any $\gamma > 0$. For any $\varphi
\in W_p^+$, let
$$\tau_+(\varphi):=\sup \{\gamma > 0 \mbox{ such that } (\mathcal{T}+\mathcal{K}(\gamma))\,\varphi
\mbox{ is nonnegative}\,\}$$ with the convention $\sup \varnothing
=0$. Then
$$\ke=\underset{\varphi \in W_p^+}{\sup}
\tau_+(\varphi).$$
\end{proposition}
%
%
\begin{proof} Let $\varphi \in X_p$ be a nonnegative eigenfunction of
$(0-\mathcal{T})^{-1}\mathcal{K}(\ke)$ associated with the spectral
radius $r_{\sigma}[(0-\mathcal{T})^{-1}\mathcal{K}(\ke)]=1,$ i.e.
$$(0-\mathcal{T})^{-1}\mathcal{K}(\ke)\,\varphi=\varphi.$$
Clearly, $\varphi \in \mathscr{D}(\mathcal{T})$. Moreover, since
$(0-\mathcal{T})^{-1}\mathcal{K}(\ke)$ is power-compact and
irreducible, a well-known consequence of  {K}rein-Rutman Theorem is
that
$$\varphi(\omega) > 0\:\:\mbox{ a.e. } \omega \in \mathbf{\Omega}$$
i.e. $\varphi \in W_p^+$. Now
$$(\mathcal{T}+\mathcal{K}(\ke))\,\varphi=0$$ is nonnegative
so that $\tau_+(\varphi) \geq \ke$ and consequently
$$\ke \leq \underset{\varphi \in W_p^+}{\sup}
\tau_+(\varphi).$$ Assume now that $\ke < \underset{\varphi \in
W_p^+}{\sup} \tau_+(\varphi)$ et let $\psi \in W_p^+$ be such that
$\tau_+(\psi) > \ke.$ Denote $\gamma=\tau_+(\psi).$ By definition
$(\mathcal{T}+\mathcal{K}(\gamma))\,\psi  \geq 0$, i.e.
$(0-\mathcal{T})^{-1}\mathcal{K}(\gamma)\,\psi \geq \psi.$ Thus, for
any $n \in \mathbb{N}$,
$[(0-\mathcal{T})^{-1}\mathcal{K}(\gamma)]^n\,\psi  \geq \psi$  so
$$r_{\sigma}[(0-\mathcal{T})^{-1}\mathcal{K}(\gamma)] \geq 1.$$
Since, $r_{\sigma}[(0-\mathcal{T})^{-1}\mathcal{K}(\ke)]=1$, Remark
\ref{rem1} ensures that $\gamma \leq \ke,$ which contradicts the
choice of $\tau_+(\psi)
> \ke.$
\end{proof}

%
The following illustrates the fact that the extremal value to the
above variational result is reached only by the nonnegative solution
to the spectral problem \eqref{abs}.
\begin{corollary}
\label{tau+1} Assume that $(0-\mathcal{T})^{-1}\mathcal{K}(\gamma)$
is power-compact and irreducible for any $\gamma > 0$. Then, for any
nonnegative $\varphi \in \mathscr{D}(\mathcal{T}) \setminus \{0\}$
\begin{equation*}
\tau_+(\varphi)=\ke \qquad \mbox{ if and only if } \qquad
(\mathcal{T}+\mathcal{K}(\ke))\,\varphi=0.
\end{equation*}
\end{corollary}
\begin{proof} Let $\varphi \in \mathscr{D}(\mathcal{T}) \setminus \{0\},$ $\varphi \geq 0,$ be
such that $\tau_+(\varphi)=\ke.$ Then
$(\mathcal{T}+\mathcal{K}(\ke))\,\varphi \geq 0$ i.e.
\begin{equation}
\label{geq} (0-\mathcal{T})^{-1}\mathcal{K}(\ke)\,\varphi \geq
\varphi.
\end{equation}
Suppose that
\begin{equation}
\label{neq} \varphi \neq
(0-\mathcal{T})^{-1}\mathcal{K}(\ke)\,\varphi,
\end{equation}
and let $\psi^\star \in X_q$ be a {\it nonnegative} eigenfunction of
the dual operator $((0-\mathcal{T})^{-1}\mathcal{K}(\ke))^\star$
associated to
$r_{\sigma}[((0-\mathcal{T})^{-1}\mathcal{K}(\ke))^\star]=
r_{\sigma}[(0-\mathcal{T})^{-1}\mathcal{K}(\ke)]=1.$ Then, by
\eqref{geq} and \eqref{neq}
\begin{equation*}
\begin{split}
\left \langle \varphi,\psi^\star \right \rangle &< \left \langle
(0-\mathcal{T})^{-1}\mathcal{K}(\ke)\,\varphi,\psi^\star
\right \rangle \\
&= \left \langle \varphi, ((0-\mathcal{T})^{-1}\mathcal{K}(\ke))^\star\,\psi^\star \right \rangle=r_{\sigma}[((0-\mathcal{T})^{-1}\mathcal{K}(\ke))^\star] \left \langle \varphi,\psi^\star \right \rangle \\
&= \left \langle \varphi,\psi^\star \right \rangle,
\end{split}
\end{equation*}
which leads to a contradiction. Hence,
$\varphi=(0-\mathcal{T})^{-1}\mathcal{K}(\ke)\,\varphi,$ i.e.
$$(\mathcal{T}+\mathcal{K}(\ke))\varphi=0.$$
Conversely, if $\varphi$ is a nonnegative eigenfunction of
$\mathcal{T}+\mathcal{K}(\ke)$ associated to the {\it null
eigenvalue}, then $\tau_+(\varphi) \geq \ge$ and the identity
$\tau_+(\varphi)=\ge$ follows from Proposition \ref{tau}.\end{proof}
\begin{remark}\label{newremark} A careful reading of the proof here above shows that, if
$\varphi \in X_p$ is such that
$(0-\mathcal{T})^{-1}\mathcal{K}(\ke)\,\varphi \geq \varphi \geq 0$,
then $\varphi \in \mathscr{D}(\mathcal{T})$ and
$\left(\mathcal{T}+\mathcal{K}(\ke) \right)\varphi=0.$\end{remark}
Let us denote by $\phi_\mathrm{eff}$ the unique critical
eigenfunction with unit norm, i.e. $\phi_\mathrm{eff} \in W_p^+$
satisfies
$$(\mathcal{T}+\mathcal{K}(\ge))\phi_\mathrm{eff}=\phi_\mathrm{eff},\qquad
\|\phi_\mathrm{eff}\|=1.$$ Then one can prove the following
approximation resolution for the criticality eigenfunction
$\phi_\mathrm{eff}$  whose proof is inspired by \cite[Theorem
7]{mkpre}. Such a result shall be  hopefully useful for practical
numerical approximation of the critical mode $\phi_\mathrm{eff}$ of
the reactor.

\begin{theorem}\label{theo:app}
Let $(\varphi_k)_k \in \mathcal{D}(\mathcal{T}) \cap X_p^+$ be such
that $\gamma_k:=\tau_+(\varphi_k) \to \ke$. We assume $\varphi_k$ to
be normalized by \begin{equation}\label{normal}
\|\left[(0-\mathcal{T})^{-1}\mathcal{K}(\gamma_k)\right]^N\varphi_k\|=1
\qquad (k \in \mathbb{N})\end{equation} where $N$ is the integer
given by Remark \ref{unifN}. Let us assume that $1$ is a simple
eigenvalue of $(0-\mathcal{T})^{-1}\mathcal{K}(\ke)$. Moreover, in
the case $p=1$, let us assume that the dual operator
$\left[(0-\mathcal{T})^{-1}\mathcal{K}(\ke)\right]^{\star}$ admits
an eigenfunction associated to its spectral radius which is bounded
away from zero. Then,
$$\lim_{k \to \infty}\|\varphi_k-\phi_\mathrm{eff}\|=0$$
where $\phi_\mathrm{eff} \in \mathcal{D}(\mathcal{T}) \cap X_p^+$ is
the unique positive eigenfunction of
$(0-\mathcal{T})^{-1}\mathcal{K}(\ke)$ associated to 1 and with unit
norm .\end{theorem}
\begin{proof} 
According to the definition of $\gamma_k:=\tau_+(\varphi_k)$,
$(\mathcal{T}+\mathcal{K}(\gamma_k))\varphi_k$ is nonnegative.
Therefore,
\begin{equation}\label{minorephik}
\varphi_k \leq (0-\mathcal{T})^{-1}\mathcal{K}(\gamma_k)\varphi_k
\end{equation}
and, iterating up to $N$ \begin{equation}\label{majoraN} \varphi_k
\leq
\left[(0-\mathcal{T})^{-1}\mathcal{K}(\gamma_k)\right]^N\varphi_k
\qquad (k \in \mathbb{N}).\end{equation} This shows, according to
\eqref{unifN}, that $\|\varphi_k\| \leq 1$. Now, if $1< p < \infty$
there exists a subsequence $(\psi_k)_k$ which converges weakly to
some $\psi \in X_p$. If $p=1$, the fact that $\gamma_k \to \ke$
combined with the compactness of
$\left[(0-\mathcal{T})^{-1}\mathcal{K}(\gamma_k)\right]^N$ lead to
the relative compactness of
$$\left\{\left[(0-\mathcal{T})^{-1}\mathcal{K}(\gamma_k)\right]^N\varphi_k\right\}_k$$ in
$X_1$. In particular, this sequence is equi--integrable and by
domination \eqref{majoraN}, $(\varphi_k)_k$ is also
equi--integrable. We can extract a subsequence $(\psi_k)_k$
converging weakly to some $\psi \in X_1$. In both cases, the
compactness  of
$\left[(0-\mathcal{T})^{-1}\mathcal{K}(\gamma_k)\right]^N$ together
with $\gamma_k \to \ke$ yield to the strong convergence
$$\left[(0-\mathcal{T})^{-1}\mathcal{K}(\gamma_k)\right]^N\psi_k \to
\left[(0-\mathcal{T})^{-1}\mathcal{K}(\ke)\right]^N\psi$$ so that
$$\|\left[(0-\mathcal{T})^{-1}\mathcal{K}(\ke)\right]^N\psi\|=1.$$
In particular $\psi \neq 0$, and, taking the weak limit in
\eqref{minorephik}, $\psi \leq
(0-\mathcal{T})^{-1}\mathcal{K}(\ke)\psi$. Now, according to Remark
\ref{newremark}, this last inequality is actually an equality, i. e.
$$\psi=(0-\mathcal{T})^{-1}\mathcal{K}(\ke)\psi.$$
Iterating again, one gets
$$\|\psi\|=\|\left[(0-\mathcal{T})^{-1}\mathcal{K}(\ke)\right]^N\psi\|=1.$$
Now, since $1$ is a simple eigenvalue of
$(0-\mathcal{T})^{-1}\mathcal{K}(\ke)$, the set of eigenfunctions of
$(0-\mathcal{T})^{-1}\mathcal{K}(\ke)$ with unit norm reduces to a
singleton. This shows that $\psi$ is the (unique) weak limit of any
subsequence of $(\varphi_k)_k$ so that the whole sequence
$(\varphi_k)_k$ {\it converges weakly} to $\psi \in X_p$. The
remainder of the proof consists in showing that the convergence
actually holds in the {\it
  strong sense}.

Let us consider first the case $1 < p < \infty$. To show now that
$\|\varphi_k - \psi\| \to 0$, it suffices to prove that
$\|\varphi_k\| \to \|\psi\|.$ A consequence of the weak convergence
leads to
$$1=\|\psi\| \leq \liminf_{k \to \infty}\|\varphi_k\|.$$
Since $\|\varphi_k\| \leq 1$ for any $k \in \mathbb{N}$, this proves
the Theorem for $1< p< \infty.$

Let us now assume $p=1$. Then $\varphi_k \to \psi$ strongly in $X_1$
if and only if the convergence holds in measure, i. e., for any $\Xi
\subset \mathbf{\Omega}$ with finite $\d\nu$--measure and every
$\epsilon > 0$
$$\lim_{k \to \infty} \d\nu\{\omega \in
\Xi\,;\,|\psi_k(\omega)-\psi(\omega)| > \epsilon \}=0.$$ Arguing by
contradiction, assume there exist $\Xi \subset \mathbf{\Omega}$ with
finite $\d\nu$--measure, a subsequence still denoted $(\psi_{k})_k$
and some $\delta > 0$ and some $\epsilon_0
> 0$ such that
\begin{equation}\label{mesure}
\d\nu\{\omega \in \Xi\,;\,|\psi_k(\omega)-\psi(\omega)| >
\epsilon_0\} \geq \delta \qquad \text{ for all } m \in
\mathbb{N}.\end{equation} Setting
$$\overline{\psi}_k=\left[(0-\mathcal{T})^{-1}\mathcal{K}(\gamma_k)\right]^N\psi_k$$
one has
$$\lim_{k \to \infty}\|\overline{\psi}_k-\psi\|=0$$
since
$\psi=\left[(0-\mathcal{T})^{-1}\mathcal{K}(\gamma_k)\right]^N\psi$.
Consequently,
$$\lim_{k \to \infty} \d\nu\{\omega \in
\Xi\,;\,|\overline{\psi}_k(\omega)-\psi(\omega)|
> \epsilon_0/2\}=0$$
and, one deduces immediately from \eqref{mesure} that
\begin{equation}\label{mesure1}
\d\nu\{\omega \in \Xi\,;\,|\psi_k(\omega)-\psi(\omega)| >
\epsilon_0/2\} \geq \delta/2 \qquad \text{ for large }
k.\end{equation} Now, let $\psi^{\star}\in
L^{\infty}(\mathbf{\Omega})$ be a positive eigenfunction of
$\left[(0-\mathcal{T})^{-1}\mathcal{K}(\ke)\right]^{\star}$
associated to the eigenvalue $1$, with $\psi^{\star}$ bounded away
from zero. One has
\begin{equation*}\begin{split} \left \langle
\overline{\psi}_k-\psi_k,\psi^{\star} \right \rangle &\geq
\int_{\{\omega \in \Xi\,;\,|\overline{\psi}_k-\psi_k|>
\epsilon_0/2\}}
\left(\overline{\psi}_k(\omega)-\psi_k(\omega)\right)\psi^{\star}(\omega)\, \d\nu(\omega)\\
&\geq \inf_{w}\psi^{\star}(\omega) \times \epsilon_0/2 \times
\delta/2=\eta >0\end{split}
\end{equation*}
or else
$$\left \langle
\left[(0-\mathcal{T})^{-1}\mathcal{K}(\gamma_k)\right]^N\psi_k,\psi^{\star}
\right \rangle - \left \langle \psi_k,\psi^{\star} \right \rangle
\geq \eta \qquad \text{ for large } k.$$ Equivalently
$$\left \langle
\psi_k,
\{[(0-\mathcal{T})^{-1}\mathcal{K}(\gamma_k)]^{\star}\}^N\psi^{\star}
\right \rangle - \left \langle \psi_k,\psi^{\star} \right \rangle
\geq \eta \qquad \text{ for large } k.$$ But,
$\{[(0-\mathcal{T})^{-1}\mathcal{K}(\gamma_k)]^{\star}\}^N\psi^{\star}$
converges strongly to
$\{[(0-\mathcal{T})^{-1}\mathcal{K}(\ke)]^{\star}\}^N\psi^{\star}$
in $L^{\infty}(\mathbf{\Omega})$, which, combined with the weak
convergence $\psi_k \rightharpoonup \psi$, implies
$$\left \langle
\psi_k,
\{[(0-\mathcal{T})^{-1}\mathcal{K}(\gamma_k)]^{\star}\}^N\psi^{\star}
\right \rangle \to \left \langle
\psi,\{[(0-\mathcal{T})^{-1}\mathcal{K}(\ke)]^{\star}\}^N\psi^{\star}
\right \rangle=\left \langle \psi,\psi^{\star}\right \rangle.$$ Now,
the contradiction follows from the fact that $\left \langle
\psi_k,\psi^{\star}\right \rangle \to \left \langle
\psi,\psi^{\star}\right \rangle.$
\end{proof}

{T}he following characterizes $\ke$ in terms of {\it sub-solution}
to the spectral problem \eqref{abs}.
\begin{proposition}
\label{tau^} Assume that $(0-\mathcal{T})^{-1}\mathcal{K}(\gamma)$
is power-compact and irreducible for any $\gamma > 0$. For any
$\varphi \in W_p^+,$ define
$$\tau_-(\varphi):=\inf \{\gamma > 0\,;\,-(\mathcal{T}+\mathcal{K}(\gamma))\varphi \in X_p^+\}$$
with the convention $\inf \varnothing = +\infty.$ Then
$$\ke=\underset{ \varphi \in W_p^+}{\inf} \tau_-(\varphi).$$
\end{proposition}
%
%
%
\begin{proof} Let $\varphi \in W_p^+$ be such that $\tau_-(\varphi) < +\infty,$ and
let $\gamma > \tau_-(\varphi).$ Then
\begin{equation}\label{<0} -(\mathcal{T}+\mathcal{K}(\gamma))\,\varphi \in X_p^+.
\end{equation}
Since $(0-\mathcal{T})^{-1}(X_p^+) \subset X_p^+$ (see Eq.
\eqref{assumptionT}), one gets that
\begin{equation}
\label{<01} \varphi
-(0-\mathcal{T})^{-1}\mathcal{K}(\gamma)\,\varphi \in W_p^+.
\end{equation}
Now, let $\psi ^\star \in X_q$ be a  nonnegative eigenfunction of
the dual operator $((0-\mathcal{T})^{-1}\mathcal{K}(\gamma))^\star$
associated with the spectral radius
$r_{\sigma}[((0-\mathcal{T})^{-1}\mathcal{K}(\gamma))^\star]=r_{\sigma}[(0-\mathcal{T})^{-1}\mathcal{K}(\gamma)]
.$ Note that $r_{\sigma}[(0-\mathcal{T})^{-1}\mathcal{K}(\gamma)] >
0$ according to Theorem \ref{bendp}. From  Krein-Rutman theorem,
$\psi^\star(\omega)
> 0$ $\d\nu$ - a.e. $\omega \in \mathbf{\Omega}.$ Thus, by (\ref{<01})
$$\left \langle \varphi-(0-\mathcal{T})^{-1}\mathcal{K}(\gamma)\,\varphi,\psi^\star \right \rangle > 0$$
i.e.
\begin{equation*}
\begin{split}
\left \langle \varphi,\psi^\star \right \rangle &> \left \langle
(0-\mathcal{T})^{-1}\mathcal{K}(\gamma)
\,\varphi,\psi^\star \right \rangle\\
&=\left \langle \varphi,
((0-\mathcal{T})^{-1}\mathcal{K}(\gamma))^\star\,\psi^\star \right
\rangle=r_{\sigma}[(0-\mathcal{T})^{-1}\mathcal{K}(\gamma)]  \left
\langle \varphi,\psi^\star \right \rangle.
\end{split}
\end{equation*}
Since $\left \langle \varphi,\psi^\star \right \rangle \neq 0,$ we
get $r_{\sigma}[(0-\mathcal{T})^{-1}\mathcal{K}(\gamma)] < 1.$ By
Remark \ref{rem1}, this means that $\gamma > \ke.$ Since $\gamma >
\tau_-(\varphi)$ is arbitrary, we obtain
$$\ke \leq \underset{\varphi \in W_p^+}{\inf} \tau_-(\varphi).$$
Conversely, let $\gamma > \ke$. Then,
$r_{\sigma}[(0-\mathcal{T})^{-1}\mathcal{K}(\gamma)] < 1,$ and
\begin{equation*}
\begin{split}
[I-(0-\mathcal{T})^{-1}\mathcal{K}(\gamma)]^{-1} &=\sum_{k=0}^{\infty} [(0-\mathcal{T})^{-1}\mathcal{K}(\gamma)]^k\\
&\geq [(0-\mathcal{T})^{-1}\mathcal{K}(\gamma)]^n\:\:\:\:\:\forall
\,n \in \mathbb{N}.
\end{split}
\end{equation*}
Given $\psi \in X_p^+,$ set
$\tilde{\psi}=(0-\mathcal{T})^{-1}\,\psi.$ Then, $\tilde{\psi} \in
W_p^+.$ Define
\begin{equation}
\label{phi}
\varphi=[I-(0-\mathcal{T})^{-1}\mathcal{K}(\gamma)]^{-1}\,\tilde{\psi}.
\end{equation}
Clearly, $\varphi \in D(\mathcal{T})$ and
$$\varphi \geq [(0-\mathcal{T})^{-1}\mathcal{K}(\gamma)]^n\,\tilde{\psi}\:\:\:\:(n \in \mathbb{N}).$$
Let $\varphi^\star \in X_q \setminus \{0\}$ be {\it nonnegative}.
Since $(0-\mathcal{T})^{-1}\mathcal{K}(\gamma)$ is {\it
irreducible}, there exists $n \in \mathbb{N}$ such that
$$\left \langle \varphi, \varphi^\star \right \rangle \geq  \left \langle
[(0-\mathcal{T})^{-1}\mathcal{K}(\gamma)]^n\,\tilde{\psi},
\varphi^\star \right \rangle > 0$$ Hence, $\left \langle \varphi,
\varphi^\star \right \rangle > 0 $ for any nonnegative  $
\varphi^\star \in X_q \setminus \{0\},$  so that $\varphi \in
X_p^+.$  Moreover, by (\ref{phi}),
$$-(\mathcal{T}+\mathcal{K}(\gamma))\,\varphi=(0-\mathcal{T})\,\tilde{\psi}=\psi \in X_p^+.$$
Hence, $\gamma \geq \tau_{-}(\varphi)$ which proves that $\ke \geq
\underset{\varphi \in W_p^+}{\inf} \tau_-(\varphi).$
\end{proof}

%
%
The following result shows that only the solution of the criticality
problem \eqref{abs} realizes the above variational characterization.
\begin{corollary}
\label{tau-1} Assume that $(0-\mathcal{T})^{-1}\mathcal{K}(\gamma)$
is power-compact and irreducible for any $\gamma > 0$. Then, for any
$\varphi \in W_p^+$
\begin{equation*}
\tau_-(\varphi)=\ke \quad \text{ if and only if } \quad
(\mathcal{T}+\mathcal{K}(\ke))\,\varphi=0.
\end{equation*}
\end{corollary}
\begin{proof} Let $\varphi \in W_p^+$ be such that $\tau_-(\varphi)=\ke.$
Then, $-(\mathcal{T}+\mathcal{K}(\ke))\,\varphi \geq 0$ so
$$ \varphi \leq (0-\mathcal{T})^{-1}\mathcal{K}(\ke) \varphi.$$
Arguing as in the proof of Corollary \ref{tau+1}, one can prove that
$\varphi=(0-\mathcal{T})^{-1}\mathcal{K}(\gamma)\,\varphi$, which
means that
$$(\mathcal{T}+\mathcal{K}(\ke))\,\varphi=0.$$
The sufficient condition follows directly from Proposition
\ref{tau^}.
\end{proof}

We are now able to characterize the effective multiplication factor
$\ke$ by means of Inf-Sup and Sup-Inf criteria, where we recall that
$\mathcal{K}=\mathcal{K}_{s}+\mathcal{K}_{f}.$
%
%
%
%
\begin{theorem}
\label{supinfabs} Under the assumptions of Theorem \ref{exist}, if
$\mathcal{K}_f(X_p^+) \subset X_p^+$ then  the criticality
eigenvalue $\ke$ is characterized by the following:
\begin{equation*}
\frac{1}{\ke}=\underset{\varphi \in W_p^+}{\min}\; \underset{\omega
\in \mathbf{\Omega}}{\esup}\,
\frac{-(\mathcal{T}+\mathcal{K}_{s})\varphi(\omega)}
{\mathcal{K}_{f}\varphi(\omega)}
=\underset{\varphi \in W_p^+}{\max}\; \underset{\omega \in
\mathbf{\Omega}}{\esinf}\,
\frac{-(\mathcal{T}+\mathcal{K}_{s})\varphi(\omega)}{\mathcal{K}_{f}\varphi(\omega)}.
\end{equation*}
\end{theorem}

%
%
%
\begin{proof}  Let $\varphi \in W_p^+$ be given,
\begin{equation*}
\begin{split}
\tau_+(\varphi)&=\sup \{ \gamma > 0\,;\,
(\mathcal{T}+\mathcal{K}_{s})\,\varphi +
\frac{1}{\gamma}\mathcal{K}_{f}\,\varphi \geq 0 \}=\sup \{ \gamma > 0\,;\,-(\mathcal{T}+\mathcal{K}_{s})\,\varphi \leq \frac{1}{\gamma}\mathcal{K}_{f}\,\varphi \}\\
&=\sup \{ \gamma >
0\,;\,-(\mathcal{T}+\mathcal{K}_{s})\,\varphi(\omega) \leq
\frac{1}{\gamma}\mathcal{K}_{f}\,\varphi(\omega)
\:\:\d\nu-\mbox{a.e. } \omega
 \in \mathbf{\Omega} \}.
\end{split}
\end{equation*}
Since $\mathcal{K}_{f}(X_p^+) \subset X_p^+$, one gets
$\displaystyle \frac{1}{\tau_+(\varphi)}= \underset{\omega \in
\mathbf{\Omega}}{ \esup}
\frac{-(\mathcal{T}+\mathcal{K}_{s})\,\varphi(\omega)}
{\mathcal{K}_{f}\,\varphi(\omega)}.$ By Proposition \ref{tau},
$$\frac{1}{\ke}=\underset{\varphi \in W_p^+}{\inf}
\underset{\omega \in \mathbf{\Omega} }{\esup}
\frac{-(\mathcal{T}+\mathcal{K}_{s})\,\varphi(\omega)}
{\mathcal{K}_{f}\,\varphi(\omega)}$$ and the infimum is attained for
the criticality eigenfunction. Similarly, let $\varphi \in W_p^+$
\begin{equation*}
\begin{split}
\tau_-(\varphi)&=\inf \{ \gamma > 0\,;\,
-(\mathcal{T}+\mathcal{K}_{s})\,\varphi -
\frac{1}{\gamma}\mathcal{K}_{f}\,\varphi \geq 0 \}=\inf \{ \gamma > 0\,;\,-(\mathcal{T}+\mathcal{K}_{s})\,\varphi \geq \frac{1}{\gamma}\mathcal{K}_{f}\,\varphi \}\\
&=\inf \{ \gamma >
0\,;\,-(\mathcal{T}+\mathcal{K}_{s})\,\varphi(\omega) \geq
\frac{1}{\gamma}\mathcal{K}_{f}\,\varphi(\omega)
\:\:\d\nu-\mbox{a.e. } \omega \in \mathbf{\Omega} \}.
\end{split}
\end{equation*}
So
$$\tau_-(\varphi)=\inf \{ \gamma > 0\,;\,\frac{-(\mathcal{T}+\mathcal{K}_{s})\,\varphi(\omega)}
{\mathcal{K}_{f}\,\varphi(\omega)} \geq \frac{1}{\gamma}
\:\:\d\nu-\mbox{a.e. } \omega \in \mathbf{\Omega} \},$$ i.e.
$$\frac{1}{\tau_-(\varphi)}=
\underset{\omega \in \mathbf{\Omega} }{\esinf}
\frac{-(\mathcal{T}+\mathcal{K}_{s})\,\varphi(\omega)}
{\mathcal{K}_{f}\,\varphi(\omega)}.$$ Using Proposition \ref{tau^},
one proves that
\begin{equation*}
\frac{1}{\ke}=\underset{\varphi \in W_p^+}{\sup}
\frac{1}{\tau_-(\varphi)}=\underset{\varphi \in W_p^+}{\sup}
\underset{\omega \in \mathbf{\Omega} }{\esinf}
\frac{-(\mathcal{T}+\mathcal{K}_{s})\,\varphi(\omega)}
{\mathcal{K}_{f}\,\varphi(\omega)},
\end{equation*}
which ends the proof.
\end{proof}

\subsection{The class of regular collision
operators}\label{sec:regular} We end this section by recalling the
class of regular collision operators introduced in kinetic theory by
M. Mokhtar-Kharroubi \cite{mmk}. This class of operators will also
be useful to study diffusion problems of type \eqref{diffu1}. We
assume here that the measure space $(\mathbf{\Omega},\d\nu)$ writes
as follows:
$$\mathbf{\Omega}=\mathcal{D} \times V, \qquad \d\nu(\omega)=\d x \otimes \d\mu(v),
\qquad \omega=(x,v) \in \mathbf{\Omega}$$ where $\d\mu$ is a
suitable Radon measure over $V$. Let $\mathcal{K} \in
\mathscr{B}(L^p(\mathbf{\Omega},\d\nu))$ be given by
\begin{equation}\label{kp}
\mathcal{K}\,:\,\varphi \longmapsto
\mathcal{K}\varphi(x,v)=\int_{V}k(x,v,v')\varphi(x,v')\d\mu(v') \in
L^p(\mathbf{\Omega},\d\nu )\end{equation} where the kernel
$k(\cdot,\cdot,\cdot)$ is measurable. For almost every $x \in
\mathcal{D}$, define
$$\widetilde{\mathcal{K}}(x) \;:\;\psi \in L^p(V,\d\mu) \longmapsto \int_Vk(x,v,v')\psi(v')\d\mu(v') \in L^p(V,\d\mu)$$
and assume that the mapping $\widetilde{\mathcal{K}}\;:\;x \in
\mathcal{D} \mapsto \widetilde{\mathcal{K}}(x) \in
\mathscr{B}(L^p(V,\d\mu))$ is strongly measurable and bounded, i.e.
$$\operatornamewithlimits{ess\,sup}_{x \in \mathcal{D}} \|\widetilde{\mathcal{K}}(x)\|_{\mathscr{B}(L^p(V,\d\mu))} < \infty.$$
The class of regular operators in $L^p$ spaces with $1 < p< \infty$
is given by the following (see \cite[Definition 4.1]{mmk}).
\begin{definition}[\textbf{Regular operator}]
Let $1 < p < \infty.$ The operator $\mathcal{K}$ defined by
\eqref{kp} is said to be regular if :
\begin{enumerate}
\item For almost every $
x \in \mathcal{D}$, $\widetilde{\mathcal{K}}(x) \in
\mathscr{B}(L^p(V,\d\mu))$ is a compact operator,
\item $\{\widetilde{\mathcal{K}}(x)\,;\,x \in \mathcal{D}\}$ is relatively compact in
$\mathscr{B}(L^p(V,\d\mu))$.
\end{enumerate} \end{definition}
In $L^1$-spaces, the definition differs a bit. We have the following
\cite{m2as}
\begin{definition}\label{defintionregul}
Let $\mathcal{K}$ be defined by \eqref{kp}. Then, $\mathcal{K}$ is
said to be a regular operator whenever
$\{\left|k(x,\cdot,v')\right|,\,(x,v') \in \mathcal{D} \times V\}$
is a relatively weakly compact subset of $L^1(V,\d\mu)$.
\end{definition}

The main interest of that classes of operators relies to the
following (see Ref. \cite{mmk} for $1 < p < \infty$ and Ref.
\cite{m2as} for a similar result whenever $p=1$):
\begin{proposition}[{\bf Approximation of regular operators}]\label{approx}   Let $1 < p < \infty$ and let $\mathcal{K}$ defined by \eqref{kp} be a regular operator
in $L^p(\mathcal{D} \times V,\d x\otimes\d\mu(v))$. Then,
$\mathcal{K}$ can be approximated in the norm operator by operators
of the form:
\begin{equation*}
\varphi \longmapsto \sum_{i \in I}\alpha_i(x)\beta_i(v)\int_V
\theta_i(v')\varphi(x,v')\d \mu(v')
\end{equation*}
where $I$ is finite, $\alpha_i \in L^{\infty}(\mathcal{D})$,
$\beta_i \in L^p(V,\d \mu)$ and $\theta_i \in L^q(V,\d \mu)$,
$1/p+1/q=1.$
\end{proposition}

%
\section{The critical transport problem}
\label{transport}

\subsection{Variational characterization}\label{VIIcatran}

This section is devoted to the determination of the effective
multiplication factor associated to the transport operator. We adopt
the notations of Section \ref{sec:regular}, namely
$\mathbf{\Omega}=\mathcal{D} \times V$ and $\d\nu(x,v)=\d x \otimes
\d\mu(v)$. Throughout this section, we assume $\mathcal{D}$ to be a
\textbf{\textit{convex and bounded}} open subset of $\mathbb{R}^N$
while
 $\mu$ is the Lebesgue measure over
$\mathbb{R}^N$ or on spheres.  In particular, our results cover
\textbf{\textit{continuous or multi-group neutron transport
problems}} but do not apply to transport problems with discrete
velocities. Let
$$\Gamma_-:=\bigg\{(x,v) \in \partial \mathcal{D} \times V\,;\, v \cdot n(x) < 0\bigg\}$$
where $n(x)$ denotes the outward unit normal at $x \in \partial \,
\mathcal{D}.$ Let $\mathcal{T}$ be the unbounded absorption operator
\begin{equation*}
\begin{cases}
\mathcal{T}\::\: \mathscr{D}(\mathcal{T}) \subset &X_p \longrightarrow X_p\\
&\varphi \longmapsto \mathcal{T}\varphi(x,v):=-v \cdot \nabla_x
\varphi(x,v)- \sigma(x,v)\varphi(x,v),
\end{cases}
\end{equation*}
with domain
$$\mathscr{D}(\mathcal{T})=\bigg\{\psi \in X_p\,;\:\:v \cdot \nabla_x \psi \in
X_p\: \mbox{ and }\: \psi_{|\Gamma_-}=0\bigg\}.$$ Here, the {\it
nonnegative} function $\sigma(\cdot,\cdot) \in L^{\infty}
(\mathcal{D} \times V)$ is
the \textit{collision frequency}. 
It is assumed to admit a {\it positive lower bound}
\begin{equation}
\label{sigmino} \sigma(x,v) \geq c > 0 \quad \mbox{  a.e. } (x,v)
\in \mathcal{D} \times V.
\end{equation}
Define the (full) collision operator $\mathcal{K}$ as the
\textit{bounded} linear (partial) integral operator
$$\mathcal{K} \::\: \psi \in X_p \longmapsto \mathcal{K}\psi(x,v):=\int_{V}\Sigma(x,v,v')\psi(x,v')\d\mu(v') \in X_p.$$
The collision kernel $\Sigma(\cdot,\cdot,\cdot)$ is assumed to be
nonnegative. In nuclear reactor theory, in a fissile material, this
collision kernel splits as
\begin{equation*} \Sigma(x,v,v')=\Sigma_s(x,v,v')+\Sigma_f(x,v,v')
\end{equation*}
where $\Sigma_s(x,v,v')$ describes the pure scattering phenomena and
$\Sigma_f(x,v,v')$ describes the fission processes. Define the
corresponding linear operators
$$\mathcal{K}_{s} \::\: \psi \in X_p \longmapsto \mathcal{K}_{s}\psi(x,v):=
\int_{V}\Sigma_s(x,v,v')\psi(x,v')\d\mu(v') \in X_p$$ and
$$\mathcal{K}_{f} \::\: \psi \in X_p \longmapsto \mathcal{K}_{f}\psi(x,v):=
\int_{V}\Sigma_f(x,v,v')\psi(x,v')\d\mu(v') \in X_p.$$
As we told it in Introduction, we are interested here in the
critical problem:
\begin{multline}\label{spectpb}
v \cdot \nabla_x \varphi(x,v) +\sigma(x,v)\varphi(x,v)-\int_V\Sigma_s(x,v,v')\varphi(x,v')\d\mu(v')\\
=\frac{1}{\ge}\,\int_V\Sigma_f(x,v,v')\varphi(x,v')\d\mu(v'),\end{multline}
where the eigenfunction $\varphi$ is \textit{nonnegative} and
satisfies the boundary condition $\varphi_{|\Gamma_-}=0.$  We recall
that the spectral bound of $\mathcal{T}$ is given by \cite{voigt}
$$s(\mathcal{T})=-\lim_{t \to \infty} \underset{ t < \tau(x,v)}{
\underset{(x,v) \in \mathcal{D} \times V}{\inf}}\;
t^{-1}\int_0^t\sigma(x+sv,v)\d s,
$$
with $\tau(x,v):=\inf \{s>0\;;\;x-sv \notin \mathcal{D}\}$.
Therefore, by \eqref{sigmino}, we have $s(\mathcal{T}) < 0.$
Moreover,
$$(0-\mathcal{T})^{-1}\varphi(x,v)=\int_0^{\tau(x,v)}\exp\left\{-\int_0^t\sigma(x-vs,v)ds\right\}\varphi(x-vt,v)dt,
$$
so that $(0-\mathcal{T})^{-1}$ fulfills \eqref{assumptionT}. Let us
now recall the irreducibility properties of
$(0-\mathcal{T})^{-1}\mathcal{K}_{f}$ for the continuous and
multigroup models. The following result may be found in Ref.
\cite{mmk}, Theorem 5.15, Theorem 5.16, (see also \cite{voigt}).
\begin{theorem}
\label{irredu} Let $\mathcal{D}$ be convex. Then,
$(0-\mathcal{T})^{-1}\mathcal{K}_{f}$ is irreducible in the two
following cases: \begin{enumerate}  \item $V$ is a closed subset of
$\mathbb{R}^N$ equipped with the  Lebesgue measure $\d\mu$ and there
exist $0 < c_1 < c_2 < \infty,$ such that $V_0 =\{v \in
\mathbb{R}^N\,;\,c_1 < |v| < c_2\} \subset V$ with
\begin{equation}\label{kf}
\Sigma_f(x,v,v') > 0 \quad \text{ a.e. } (x,v,v') \in
\big(\mathcal{D} \times V \times V_0\big) \cup \big(\mathcal{D}
\times V_0 \times V\big).\end{equation}
\item $V$ is the union of $k$ disjoint spheres $(k \geq 1)$,
$$V=\bigcup_{i=1}^k V_i, \qquad V_i=\{v \in \mathbb{R}^N\,;|v|=r_i\}, \qquad (r_i >0,\:i=1,\ldots,k)$$
and, on each sphere, $\d\mu$ is the surface Lebesgue measure.
Moreover, for any $i, j \in \{1,\ldots,k\},$ there exists $\ell \in
\{1,\ldots,m\}$ such that
\begin{equation}\label{kfmult}
\Sigma_f(x,v,v')
> 0 \qquad \text{ a.e. } (x,v,v') \in \big(\mathcal{D} \times V_i
\times V_\ell\big) \cup \big(\mathcal{D} \times V_\ell \times
V_j\big).\end{equation}
\end{enumerate}\end{theorem}
\begin{remark} In the above case \textit{(1)}, corresponding to \textbf{\textit{continuous models}}, it is possible to
provide different criteria ensuring the irreducibility of
$(0-\mathcal{T})^{-1}\mathcal{K}_f$ (see for instance Ref.
\cite{greiner}). In the second case \textit{(2)}, which corresponds
to \textbf{\textit{multigroup}} transport equation, several
different criteria also exist \cite{nagel}.
\end{remark}
Using the notations of Section \ref{caraabs}, we have the following
characterization of the effective multiplication factor of the
transport operator. 
%
\begin{theorem}\label{princitran}
Let us assume that $\mathcal{K}$ is a regular collision operator and
that one of the hypothesis of Theorem \ref{irredu} holds. The
critical problem \eqref{spectpb} admits a unique solution $\ge$ if
and only if
\begin{equation*} \lim_{\gamma \to
0}r_{\sigma}[(0-\mathcal{T})^{-1}\mathcal{K}(\gamma)] >1 \qquad
\mbox{ and } \qquad r_{\sigma}[(0-\mathcal{T})^{-1}\mathcal{K}_{s}]
< 1.
\end{equation*} Moreover,
 \begin{equation}%
\label{carac}
\begin{split}
\frac{1}{\ge}&=\underset{\varphi \in W_p^+}{\min}\; \underset{(x,v)
\in \mathcal{D} \times V}{\esup} \frac{v \cdot \nabla_x
\varphi(x,v)+
\sigma(x,v)\varphi(x,v)-\displaystyle{\int_{V}}\Sigma_s(x,v,v')\varphi(x,v')\,\d\mu(v')}
{\displaystyle{\int_{V}}\Sigma_f(x,v,v')\varphi(x,v')\,\d\mu(v')}\\
&=\underset{\varphi \in W_p^+}{\max}\; \underset{(x,v) \in
\mathcal{D} \times V}{\esinf} \frac{v \cdot \nabla_x \varphi(x,v)+
\sigma(x,v)\varphi(x,v)-\displaystyle{\int_{V}}\Sigma_s(x,v,v')\varphi(x,v')\,\d\mu(v')}
{\displaystyle{\int_{V}}\Sigma_f(x,v,v')\varphi(x,v')\,\d\mu(v')}.
\end{split}
\end{equation}
\end{theorem}
\begin{proof}
Since $\mathcal{K}$ is a regular collision operator, one deduces
from \cite[Theorems 4.1 \& 4.4]{mmk} when $1 < p < \infty$
(respectively  \cite{m2as} if $p=1$) that
$(0-\mathcal{T})^{-1}\mathcal{K}$ is a power-compact operator in
$X_p$ $(1 \leq p < \infty)$ under our assumptions on the measure
$\mu$. Moreover, in the continuous case, thanks to \eqref{kf}, for
any $\varphi \in X_p^+$
\begin{equation*}\begin{split}
\mathcal{K}_{f}\varphi(x,v)&=\int_{V}\Sigma_f(x,v,v')\,\varphi(x,v')\d\mu(v') \\
 &\geq \int_{V_0}\Sigma_f(x,v,v')\,\varphi(x,v')\d\mu(v')>0\qquad \text{ a.e. } (x,v) \in \mathcal{D} \times V,\end{split}\end{equation*}
i.e. $\mathcal{K}_{f}(X_p^+) \subset X_p^+.$ Similarly, in the
multigroup case, Eq. \eqref{kfmult} implies $\mathcal{K}_f(X_p^+)
\subset X_p^+$. Now, the existence of $\ke$ follows from Theorem
\ref{exist} while \eqref{carac} follows from Theorem
\ref{supinfabs}.\end{proof}

\begin{remark} Denote by $\phi_{\mathrm{eff}}$ the
nonnegative solution of \eqref{spectpb}, one can check that
$\phi_{\text{eff}} \in W_p^+$. Therefore, in \eqref{carac}, the
supremum and the infimum are reached for
$\phi=\phi_{\text{eff}}$.\end{remark}

\begin{remark}
Note that it is possible to provide  practical criteria that are
satisfied in nuclear reactor theory and that ensure the existence of
$\ge$ \cite{bal,thevenot}. Such criteria usually rely on dissipative
properties of the pure scattering operator.
\end{remark}

\begin{remark}\label{boundary} It is important to point out that the
above characterization is not restricted to the case of absorbing
conditions but also holds for general boundary conditions modeled by
some suitable \textit{\textbf{nonnegative albedo operator}}.
Actually, if one considers a transport operator $\mathcal{T}_H$
associated to general nonnegative albedo boundary operator  $H$
which relates the incoming and outgoing fluxes in $\mathcal{D}$
\cite{latrach}, then the above theorem holds true provided
$(0-\mathcal{T}_H)^{-1}\mathcal{K}$ is a power-compact operator in
$X_p$ $(1 \leq p < \infty)$ when $\mathcal{K}_s$ and $\mathcal{K}_f$
are regular operators. This is always the case whenever $1 < p <
\infty$ by virtue of the velocity averaging lemma \cite{latrach}.
The problem is more delicate in a $L^1$-setting and is related to
the geometry of the domain $\mathcal{D}$ \cite{sanchezttsp}.
\end{remark}

\subsection{Necessary conditions of super-criticality and
sub-criticality}\label{lowtran}

We shall use the result of the previous section to derive necessary
conditions ensuring the reactor to be super-critical or
sub-critical. Note that, for practical implications, a nuclear
reactor can be operative and create energy only when slightly
super-critical (i.e. $1< \ge < 1+\delta$ with $\delta >0$ small
enough), in this case, the whole chain fission being controlled by
rods of absorbing matter. Throughout this section, we shall assume
$\ke$ to exist.

We shall provide lower and upper bounds on the effective
multiplicative factor $\ge$ only when the velocity space $V$ is
\textit{\textbf{bounded away from zero}}. Recall that, since $V$ is
assumed to be closed, this means that $0 \notin V$ (see also Remark
\ref{OinV}).

For almost every $x \in \mathcal{D}$, define $\mathcal{K}^{\tau}(x)$
as the following operator on $L^p(V,\d\mu)$:
$$\mathcal{K}^\tau (x)\::\psi \in L^p(V,\d\mu) \mapsto
\int_{V}\dfrac{\Sigma(x,v,v')\tau(x,v')}{1+\sigma(x,v)\tau(x,v)}\psi(v')\d\mu(v')
\in L^p(V,\d\mu)$$ where we recall that
$\Sigma(x,v,v')=\Sigma_s(x,v,v')+\Sigma_f(x,v,v')$ and $\tau(x,v)$
is the stay time in $\mathcal{D}$. Then, one defines as in
\cite{mk}, the following
$$\overline{\vartheta}:=\inf_{\psi \in L^p_+(V,\d\mu)}\esup_{(x,v) \in \mathcal{D} \times V} \dfrac{[\mathcal{K}^\tau (x)\psi](v)}{\psi(v)}$$
where $L_+^p(V,\d\mu)=\{\psi \in L^p(V,\d\mu)\,;\,\psi(v) >0
\:\d\mu-\text{a.e. } v \in V\}.$ Then, one has the following
estimate:
\begin{proposition}\label{estimate1} Under the assumptions of Theorem \ref{princitran}, if
$\overline{\vartheta} < 1$, then $\ge \leq \overline{\vartheta}.$
\end{proposition}
\begin{proof} Assume  $\overline{\vartheta} <1$. Given $ \vartheta \in (\overline{\vartheta},1)$, let $\psi_0 \in
L_+^p(V,\d\mu)$ be such that $$\esup_{(x,v) \in \mathcal{D} \times
V} \dfrac{[\mathcal{K}^\tau (x)\psi_0](v)}{\psi_0(v)} \leq
{\vartheta}.$$ Let us consider then the following test-function
$\varphi_0(x,v)=\tau(x,v)\psi_0(v)$. Since $0 \notin V$,
$\tau(\cdot,\cdot)$ is bounded and such an application $\varphi_0$
belongs to $W^+_p$ since
$$\tau(x+tv,v)=\tau(x,v)+t \qquad \qquad \text{ a.e. } (x,v) \in \mathcal{D}
\times V, \:t \geq 0,$$ implies $v \cdot \nabla_x
\varphi_0(x,v)=\psi_0(v).$ Then, for any $\gamma >0$, one sees that
\begin{multline}\label{t+kgam}
-(\mathcal{T}+\mathcal{K}(\gamma))\varphi_0(x,v)={{\vartheta}}^{\,-1}\left(1+\sigma(x,v)\tau(x,v)\right)
\bigg({\vartheta}\psi_0(v)-[\mathcal{K}^\tau (x)\psi_0](v)\\
 +(1-{\vartheta})\int_{V}\dfrac{\Sigma_s(x,v,v')\tau(x,v')}{1+\sigma(x,v)\tau(x,v)}\psi_0(v')\d\mu(v') +\\
 +\frac{\gamma-{\vartheta}}{\gamma}\int_{V}\dfrac{\Sigma_f(x,v,v')\tau(x,v')}{1+\sigma(x,v)\tau(x,v)}\psi_0(v')\d\mu(v')\bigg).\end{multline}
Since $\Sigma_s \geq 0$ and $1-{\vartheta} \geq 0$, one sees that
\begin{multline*}-(\mathcal{T}+\mathcal{K}(\gamma))\varphi_0(x,v) \geq {{\vartheta}}^{\,-1}
\left(1+\sigma(x,v)\tau(x,v)\right)\bigg({\vartheta}\psi_0(v)-[\mathcal{K}^\tau (x)\psi_0](v)\\
 +\frac{\gamma-{\vartheta}}{\gamma}\int_{V}\dfrac{\Sigma_f(x,v,v')\tau(x,v')}{1+\sigma(x,v)\tau(x,v)}\psi_0(v')\d\mu(v')\bigg).\end{multline*}
In particular, from the positivity of $\Sigma_f$, one sees that,
provided $\gamma \geq {\vartheta}$,
$-(\mathcal{T}+\mathcal{K}(\gamma))\varphi_0(x,v) \geq 0$ for almost
every $(x,v) \in \mathcal{D} \times V$. Then, from Proposition
\ref{tau^}, this means that $\tau_-(\varphi_0) \leq {\vartheta}$ and
$\ge \leq {\vartheta}$. Since $\vartheta > \overline{\vartheta}$ is
arbitrary, one gets the result.\end{proof}
\begin{remark} From the above result, one sees that the reactor is
sub-critical whenever $\overline{{\vartheta}} < 1$. Note that the
fact that $\overline{\vartheta} <1$ implies $\ge \leq 1$ is already
contained in \cite[Theorem 7]{mk}.
\end{remark}
The above result provides an upper bound of $\ge$ leading to the
sub-criticality of the reactor core. To get a lower bound of $\ge$,
one defines a similar quantity
$$\underline{\vartheta}=\sup_{\psi \in L^p_+(V,\d\mu)}\esinf_{(x,v) \in \mathcal{D} \times V} \dfrac{[\mathcal{K}^\tau
(x)\psi](v)}{\psi(v)}.$$
\begin{proposition}\label{estimate1bis} Under the assumptions of Theorem
\ref{princitran}, if $\underline{\vartheta} > 1$, then $\ge \geq
\underline{\vartheta}.$ In particular, for a reactor core to be
sub-critical, it is necessary that $\underline{\vartheta} \leq 1$.
\end{proposition}
\begin{proof} The proof is very similar to that of Prop.
\ref{estimate1}. Namely, assume $\underline{\vartheta} >1$. For any
$\vartheta \in (1,\underline{\vartheta})$, let
 $\psi_0 \in L_+^p(V,\d\mu)$ be such
that $ [\mathcal{K}^\tau (x)\psi_0](v) \geq \vartheta  {\psi_0(v)}$
for almost every $(x,v) \in \mathcal{D} \times V$. Then, the
function $\varphi_0(x,v)=\tau(x,v)\psi_0(v)$ belongs to $W_p^+$ and,
arguing as in Prop. \ref{estimate1}, thanks to Eq. \eqref{t+kgam},
one sees that, $(\mathcal{T}+\mathcal{K}(\gamma))\varphi$ is
nonnegative for any $\gamma \leq \vartheta$. Consequently,
$\tau_+(\varphi) \geq \vartheta$ and Prop. \ref{tau} implies that
$\ke \geq \vartheta$ for any $\vartheta \in
(1,\underline{\vartheta}).$
\end{proof}

\begin{remark} To the author's knowledge, the identity
$\overline{\vartheta}=\underline{\vartheta}$ is an open question.
Notice however  that,  according to I. Marek's result, Ref.
\cite{marek1}, Theorem 3.2, for any $x \in \mathcal{D}$, one has the
identity
$$\sup_{\psi \in L^p_+(V,d\mu)}\esinf_{v  \in  V} \dfrac{[\mathcal{K}^\tau
(x)\psi](v)}{\psi(v)}=\inf_{\psi \in L^p_+(V,d\mu)}\esup_{v  \in V}
\dfrac{[\mathcal{K}^\tau
(x)\psi](v)}{\psi(v)}=r_\sigma[\mathcal{K}^\tau(x)],$$ where we
recall that $\mathcal{K}^\tau(x)$ is an operator in $L^p(V,\d\mu)$.
\end{remark}
In the same spirit, for almost every $x \in \mathcal{D}$, define
$\mathcal{K}^{\tau}_f(x)$ as the following operator on
$L^p(V,\d\mu)$:
$$\mathcal{K}_{f}^{\tau} (x)\::\psi \in L^p(V,\d\mu) \mapsto
\int_{V}\dfrac{\Sigma_f(x,v,v')\tau(x,v')}{1+\sigma(x,v)\tau(x,v)}\psi(v')\d\mu(v')
\in L^p(V,\d\mu)$$ and let us define, as in \cite{mk}, the set $I_f$
of all $\beta \geq 0$ for which there exists $\psi \in
L^p_+(V,\d\mu)\setminus \{0\}$ such that
$$[\mathcal{K}_f^\tau(x)\psi](v) \geq \beta \psi(v),\quad
\text{ for almost every } \quad (x,v) \in \mathcal{D} \times V.$$
According to \cite[Lemma 4]{mk} the set $I$ is closed so that, if
one defines
$${\beta}_f:=\sup\{\beta, \beta \in I\}$$
then, there exists ${\psi_f} \in L^p_+(V,\d\mu)\setminus \{0\}$ such
that $[\mathcal{K}_f^\tau(x)\psi_f](v) \geq \beta_f  {\psi_f}(v)$
for almost every $(x,v) \in \mathcal{D} \times V.$ When the velocity
space is bounded away from $0$ then, ${\beta}_f$ provides a lower
bound for $\ge$:

\begin{proposition} Under the assumptions of Theorem \ref{princitran}, one has $\ge \geq {\beta}_f.$
\end{proposition}
\begin{proof} Set $\varphi_f(x,v)=\tau(x,v)\psi_f(v)$ where $\psi_f \in L^p_+(V,\d\mu) $
is defined here above. Arguing as in the proof of Proposition
\ref{estimate1}, one sees that, since $0 \notin V$, $\varphi_f \in
W_p^+$. Therefore, Theorem \ref{princitran} ensures that
$$\frac{1}{\ge}\leq \underset{(x,v)
\in \mathcal{D} \times V}{\esup} \frac{v \cdot \nabla_x
\varphi_f(x,v)+
\sigma(x,v)\varphi_f(x,v)-\displaystyle{\int_{V}}\Sigma_s(x,v,v')\varphi_f(x,v')\,\d\mu(v')}
{\displaystyle{\int_{V}}\Sigma_f(x,v,v')\varphi_f(x,v')\,\d\mu(v')}.$$
As in the proof of Prop. \ref{estimate1} and since $\Sigma_s $ and
$\varphi_s$ are nonnegative, one gets
$$\frac{1}{\ge}\leq \underset{(x,v)
\in \mathcal{D} \times V}{\esup} \frac{ \psi_f(v) }{
[\mathcal{K}_f^\tau(x)\psi_f](v)} \leq \frac{1}{{\beta}_f}$$ which
ends the proof.\end{proof}

\begin{remark} The above Proposition provides a lower bound of the
criticality eigenvalue $\ge$ that depends \textit{\textbf{only on
the fission collision operator $\mathcal{K}_f$}}. In particular, a
sufficient condition for the reactor to be super-critical is
${\beta}_f > 1$.
\end{remark}

\begin{proposition}
\label{0-V}Let $v_0:=\inf \{|v|\;;v \in V\}$ and let $d$ be the
diameter of $\mathcal{D}.$ Define
$$\mathbf{\Lambda}_f:={\underset{\psi \in L^p_+(V,\d\mu)}
{\sup}}\;\underset{(x,v) \in \mathcal{D} \times V} {\esinf}
\frac{1}{\psi(v)}\int_{V}\Sigma_f(x,v,v')\tau(x,v')\psi(v')\d\mu(v').$$
Then, under the assumptions of Theorem \ref{princitran},
\begin{equation}
\label{inf} \dfrac{1}{\ge} \leq \frac{1+\overline{\sigma}\,
\frac{d}{v_0} }{\mathbf{\Lambda}_f} ,
\end{equation}
where $\overline{\sigma}:=\underset{(x,v) \in \mathcal{D} \times V}
{\esup}\sigma(x,v).$ In particular, if $V$ bounded then
\begin{equation} \label{vborne} \ge \geq
\frac{1}{1+\overline{\sigma}\,\frac{d}{v_0}}\; \underset{(x,v) \in
\mathcal{D} \times V} {\esinf}\int_{V}
\Sigma_f(x,v,v')\tau(x,v')\d\mu(v').
\end{equation}
\end{proposition}
%
%
%
\begin{proof} Let us consider again test-functions of the form
$\varphi(x,v)=\tau(x,v)\psi(v)$ where $\psi \in L^p_+(V,\d\mu)$.
Then, as above, according to \eqref{carac}
 \begin{equation*}
 \begin{split}
\frac{1}{\ge} &\leq  \underset{(x,v) \in \mathcal{D} \times V}
{\esup} \frac{1+ \sigma(x,v)\tau(x,v)-\dfrac{1}{\psi(v)}
\displaystyle{\int_{V}}\Sigma_s(x,v,v')\tau(x,v')\psi(v')\d\mu(v')}
{\dfrac{1}{\psi(v)}\displaystyle{\int_{V}}
\Sigma_f(x,v,v')\tau(x,v')\psi(v')\d\mu(v')}\\
&\leq \frac{1+\underset{(x,v) \in \mathcal{D} \times V}
{\esup}\sigma(x,v)\tau(x,v)} {\underset{(x,v) \in \mathcal{D} \times
V} {\esinf}
\dfrac{1}{\psi(v)}\displaystyle{\int_{V}}\Sigma_f(x,v,v')\tau(x,v')\psi(v')\d\mu(v')}.
\end{split}\end{equation*} Since such an inequality holds for arbitrary
$\psi(v) >0$ and since $\tau(x,v) \leq d /|v|  \leq  d/ v_0$ for
almost every $(x,v) \in \mathcal{D} \times V,$  one gets
\eqref{inf}. To prove \eqref{vborne}, it suffices to consider the
test-function
 $\psi(v)=1$ $(v \in V),$ which belongs to $L^p_+(V,\d\mu)$
provided $V$ is bounded.
\end{proof}
\begin{remark}\label{OinV} We dealt in this section with the case of velocities
bounded away from zero. For practical use in nuclear engineering,
this is no major restriction. However, it should also be possible to
derive explicit bounds of $\ge$ when $0 \in V$. In such a case, the
exit time $\tau(x,v)$ is not bounded anymore but behave as
$\frac{1}{|v|}$ for small $|v|$. Therefore, test-functions of the
form $\varphi(x,v)=\tau(x,v)\psi(v)$ belong to $W_p^+$ if and only
if $\frac{\psi(v)}{|v|} \in L^p_+(V,\d\mu)$.
\end{remark}
%
%
\section{The critical problem for the energy-dependent diffusion model}\label{diffusion}

\subsection{Variational characterization}
In this section, we are concerned with the following
\begin{multline}\label{diffu}
-\di\left(\mathbb{D}(x,\xi)\nabla_x \varrho(x,\xi)\right)+\sigma(x,\xi)\varrho(x,\xi)-\int_{E}\Sigma_s(x,\xi,\xi')\,\varrho(x,\xi')\,\d\xi'\\
=\frac{1}{\ke}\int_E\Sigma_f(x,\xi,\xi')\,\varrho(x,\xi')\,\d\xi',\end{multline}
 where the unknown $\varrho(\cdot,\cdot)$ is assumed to be \textit{nonnegative} and to
 satisfy the Dirichlet boundary conditions
\begin{equation*}\label{dirichlet}
\varrho_{|\partial \mathcal{D}}(\cdot,\xi)=0\qquad \qquad \text{
a.e. } \xi \in E,\end{equation*} where $\mathcal{D}$ is $C^2$ open
\textit{bounded and connected }subset of $\mathbb{R}^N$ and $E$ is
an interval of $]0,\infty[$. We will assume throughout this section
that there exist some constants $\sigma_i>0$ $(i=1,2)$ such that
\begin{equation}\label{dsig}
0<\sigma_1 \leq \sigma(x,\xi) \leq \sigma_2 < \infty,\\
\text{ a.e. } (x,\xi) \in \mathcal{D} \times E.\end{equation}
Moreover, we assume the measurable matrix--valued application
$\mathbb{D}(\cdot,\cdot)$ satisfies the following (uniform)
ellipticity property
\begin{equation}\label{ellip}
\esinf_{(x,\xi) \in \mathcal{D}\times E
}\sum_{i,j=1}^Nd_{ij}(x,\xi)\eta_i\eta_j \geq d_1|\eta|^2 \qquad
\qquad (\eta \in \mathbb{R}^N)\end{equation} and regularity
assumption $d_{ij}(\cdot,\xi) \in W^{1,2}_\mathrm{loc}(\mathcal{D})$
for almost every $\xi \in E$. We will study Problem \eqref{diffu} in
a \textit{Hilbert space setting} for simplicity. Namely, set
$$X_2=L^2(\mathcal{D} \times E,\d x\d\xi).$$
Let us assume the kernels $\Sigma_s(\cdot,\cdot,\cdot)$ and
$\Sigma_f(\cdot,\cdot,\cdot)$ to be {\it nonnegative} and define the
\textit{scattering operator}
$$\mathcal{\mathcal{K}}_{s}\::\:\psi \in X_2 \longmapsto \mathcal{K}_{s}\psi(x,\xi)=\int_E\Sigma_s(x,\xi,\xi')\psi(x,\xi')\d\xi' \in X_2,$$
and the \textit{fission operator}
$$\mathcal{K}_{f}\::\:\psi \in X_2 \longmapsto \mathcal{K}_{f}\psi(x,\xi)=\int_E\Sigma_f(x,\xi,\xi')\psi(x,\xi')\d\xi' \in X_2.$$
We will assume $\mathcal{K}_{s}$ and $\mathcal{K}_{f}$ to be {\it
bounded} operators in $X_2.$ Define then the full collision operator
$$\mathcal{K}\::\:\psi \in X_2 \longmapsto \mathcal{K}\psi(x,\xi)=\int_E\Sigma(x,\xi,\xi')\psi(x,\xi')\d\xi' \in X_2,$$
where
$$\Sigma(x,\xi,\xi')=\Sigma_s(x,\xi,\xi')+\Sigma_f(x,\xi,\xi') \qquad(x,\xi,\xi') \in \mathcal{D} \times E \times E.$$
Let us introduce the \textit{diffusion operator}
\begin{equation*}\begin{cases}
\mathcal{T}\::\:&\mathscr{D}(\mathcal{T}) \subset X_2 \longrightarrow X_2\\
&\varrho \longmapsto
\mathcal{T}\varrho(x,\xi)=\di(\mathbb{D}(x,\xi)\nabla_x\varrho(x,\xi))-\sigma(x,\xi)\varrho(x,\xi),\end{cases}\end{equation*}
with domain
$$\mathscr{D}(\mathcal{T})=\{\psi \in X_2\,;\,\psi(\cdot,\xi) \in  {H}_0^1(\mathcal{D}) \cap  {H}^2(\mathcal{D})\:\text{ a.e. } \xi \in E \text{ and } \mathcal{T}\psi \in X_2 \}$$
where $ {H}_0^{1}(\mathcal{D})$ and $ {H}^{2}(\mathcal{D})$ are the
usual Sobolev spaces. With these notations, the spectral problem
\eqref{diffu} reads
$$(\mathcal{T}+\mathcal{K}_{s}+\frac{1}{\ke}\mathcal{K}_{f})\varrho_{\mathrm{eff}}=0, \qquad \qquad \varrho_{\mathrm{eff}} \in \mathscr{D}(\mathcal{T})\;,\,\varrho_{\mathrm{eff}} \geq 0,\:\varrho_{\mathrm{eff}} \neq 0.$$
According to the strong maximum principle, it is clear that
$s(\mathcal{T}) < 0$ and $(0-\mathcal{T})^{-1}(X_2^+) \subset
X_2^+$. In order to apply Theorem \ref{supinfabs}, one has to make
sure that $(0-\mathcal{T})^{-1}\mathcal{K}$ is power-compact and
that $(0-\mathcal{T})^{-1}\mathcal{K}_{f}$ is irreducible. Let us
begin with the following  compactness result which is similar to the
usual {\it velocity averaging} lemma (see  \cite{golse} and
\cite[Chapter 2]{mmk}) for transport equations and is based on some
consequence of the Sobolev embedding Theorem \cite{brezis}.
\begin{theorem}\label{compacite}
If $\mathcal{K} \in \mathscr{B}(X_2)$ is regular then
$\mathcal{K}(0-\mathcal{T})^{-1}$ is a compact operator in
$X_2$.\end{theorem}
\begin{proof} By Proposition \ref{approx}, it suffices to prove the result for a collision operator $\mathcal{K}$
such that
$$\mathcal{K}\::\:\varrho \in X_2 \longmapsto \mathcal{K}\varrho(x,\xi)=\alpha(x)h(\xi)\int_Ef(\xi')\varrho(x,\xi')\d\xi' \in X_2$$
where
$$\alpha \in L^{\infty}(\mathcal{D}), \quad h \in L^2(E,\d\xi)\:\text{ and } f \in L^2(E,\d\xi).$$
Moreover, by a density argument, one can also assume $f$ and $h$ to
be continuous functions with compact support in $E.$  Let us split
$\mathcal{K}(0-\mathcal{T})^{-1}$ as:
$$\mathcal{K}(0-\mathcal{T})^{-1}=\Theta \mathcal{M}(0-\mathcal{T})^{-1}$$
where
$$\Theta\;:\:\varrho \in L^2(\mathcal{D},\d x) \longmapsto [\Theta\varrho](x,\xi)=\alpha(x)h(\xi)\varrho(x) \in X_2,$$
and $\mathcal{M}$ is the \textit{averaging operator}
$$\mathcal{M}\;:\;\psi \in X_2 \longmapsto \mathcal{M}\psi(x)=\int_{E}f(\xi')\psi(x,\xi')\d\xi' \in L^2(\mathcal{D}).$$
It is enough to prove that $\mathcal{M}(0-\mathcal{T})^{-1}\,:\,X_2
\to L^2(\mathcal{D})$ is compact. Let $\mathcal{B}$ be a bounded
subset of $X_2$. One has to show that $\{\mathcal{M}g\,;\,g \in
(0-\mathcal{T})^{-1}(\mathcal{B})\}$ is a relatively compact subset
of $L^2(\mathcal{D}).$ For any $\varphi \in \mathcal{B}$, set
$$g(x,\xi)=(0-\mathcal{T})^{-1}\varphi(x,\xi).$$
For almost every $\xi \in E,$ $g(\cdot,\xi) \in H^1_0(\mathcal{D}).$
One extends $g$ to the whole space $\mathbb{R}^N$ by
\begin{equation*}
\tilde{g}(x,\xi)=\begin{cases}g(x,\xi)\:\text{if } x \in \mathcal{D}\\
0 \;\text{ else}.\end{cases}\end{equation*} Clearly, for almost
every $\xi \in E$, $\tilde{g}(\cdot,\xi) \in H^1(\mathbb{R}^N).$
Consequently, according to \cite[Proposition IX.3]{brezis}, for a.
e. $\xi \in E$ and any $h \in \mathbb{R}^N$
\begin{equation*}
\|\tau_h\tilde{g}(\cdot,\xi)-\tilde{g}(\cdot,\xi)\|_{L^2(\mathcal{D})}
\leq |h|\|\nabla_x\tilde{g}(\cdot,\xi)\|_{L^2(\mathcal{D})},
\end{equation*}
where $\tau_hf(x)=f(x+h)$ $(x \in \mathcal{D},\,h \in
\mathbb{R}^N)$, i.e.
\begin{equation}\label{tauh}
\int_{\mathcal{D}}|\tilde{g}(x+h,\xi)-\tilde{g}(x,\xi)|^2\d x \leq
|h|^2\int_{\mathcal{D}}|\nabla_x{g}(x,\xi)|^2\d x.
\end{equation}
Now, recall that
$$-\di(\mathbb{D}(x,\xi)\nabla_x g(x,\xi))+\sigma(x,\xi)g(x,\xi)=\varphi(x,\xi) \qquad (x,\xi) \in \mathcal{D} \times E.$$
Multiplying this identity by $g(x,\xi)$ and integrating by parts
yield, thanks to the ellipticity property \eqref{ellip},
\begin{equation*}
 d_1 \int_{\mathcal{D} \times E}|\nabla_xg(x,\xi)|^2\d x\d\xi \leq \int_{\mathcal{D} \times E}|g(x,\xi)||\varphi(x,\xi)|\d x\d\xi.\end{equation*}
In particular, since $\mathcal{B}$ is bounded, by Cauchy-Schwarz
inequality, there exists $c
> 0$ such that
\begin{equation}\label{borne}
\sup_{g \in (0-\mathcal{T})^{-1}(\mathcal{B})}\int_{\mathcal{D}
\times E}|\nabla_xg(x,\xi)|^2\d x\d\xi \leq c.\end{equation} Then,
\eqref{borne} together with \eqref{tauh} yield
$$\int_{\mathcal{D} \times E}|\tilde{g}(x+h,\xi)-\tilde{g}(x,\xi)|^2\d x\d\xi \leq c|h|^2.$$
By H\"older's inequality, since $f$ is continuous with compact
support
\begin{equation*}\begin{split}
\int_{\mathcal{D}}|M\tilde{g}(x+h)-M\tilde{g}(x)|^2\d
x&=\int_{\mathcal{D}}\d x
\left|\int_E(\tilde{g}(x+h,\xi)-\tilde{g}(x,\xi))f(\xi)\d\xi\right|^2\\
&\leq C \int_{\mathcal{D} \times
E}|\tilde{g}(x+h,\xi)-\tilde{g}(x,\xi)|^2\d x\d\xi\\ &\leq
|h|^2C,\end{split}\end{equation*} where $C>0$ does not depend on
$g$. In particular,
$$\lim_{h \to 0}\sup_{g \in (0-\mathcal{T})^{-1}(\mathcal{B})}\int_{\mathcal{D}} |M\tilde{g}(x+h)-M\tilde{g}(x)|^2\d x=0.$$
Now, using that $\widetilde{Mg}=M\tilde{g}$ one deduces the
conclusion from Riesz-Fr\'echet-{K}olmogorov Theorem
\cite{brezis}.\end{proof}

We are now in position to prove the main result of this section
where the notations of Section \ref{caraabs} are adopted :
\begin{theorem}\label{prindiff}
Let $\mathcal{K} \in \mathscr{B}(X_2)$ be regular. Assume there
exists an open subset $E_0 \subset E$ such that
\begin{equation}\label{kap_f}
\Sigma_f(x,\xi,\xi') > 0\:\:\:\:\text{ a.e. } (x,\xi,\xi') \in
\mathcal{D} \times E \times E_0.\end{equation} Then, the problem
\eqref{diffu} admits a effective multiplication factor $\ke
> 0$ if, and only if,
$$\lim_{\gamma \to 0}r_{\sigma}[(0-\mathcal{T})^{-1}\mathcal{K}(\gamma)] >1 \qquad \mbox{ and
} \qquad  r_{\sigma}[(0-\mathcal{T})^{-1}\mathcal{K}_{s}] < 1.$$
Moreover, $\ke$ is characterized {\small \begin{equation}
\label{caracdiff}
\begin{split}
\frac{1}{\ke}&=\underset{\varphi \in W_p^+}{\min}\; \underset{(x,v)
\in \mathcal{D} \times E}{\operatornamewithlimits{ess\,sup}}
\frac{-\mathrm{div}(\mathbb{D}(x,\xi) \nabla_x\varphi(x,\xi))+
\sigma(x,\xi)\varphi(x,\xi)-\displaystyle{\int_E}\Sigma_s(x,\xi,\xi')\varphi(x,\xi')\d\xi'}
{\displaystyle{\int_E}\Sigma_f(x,\xi,\xi')\varphi(x,\xi')\d\xi'}\\
&=\underset{\varphi \in W_p^+}{\max}\; \underset{(x,\xi) \in
\mathcal{D} \times V}{\operatornamewithlimits{ess\,inf}}
\frac{-\mathrm{div}(\mathbb{D}(x,\xi) \nabla_x \varphi(x,\xi))+
\sigma(x,\xi)\varphi(x,\xi)-\displaystyle{\int_E}\Sigma_s(x,\xi,\xi')\varphi(x,\xi')\d\xi'}
{\displaystyle{\int_{E}}\Sigma_f(x,\xi,v')\varphi(x,\xi')\d\xi'}.
\end{split}
\end{equation}}
\end{theorem}
\begin{proof} From \eqref{kap_f}, $\mathcal{K}_f(X_2^+) \subset X_2^+$. Now, for almost every $\xi \in E$, define $T_\xi$ as
the following operator on $L^2(\mathcal{D})$:
$$T_\xi\::\:\varrho \in \mathscr{D}(T_\xi) \mapsto T_\xi \varrho(x)=
\di(\mathbb{D}(x,\xi)\nabla_x\varrho(x))-\sigma(x,\xi)\varrho(x),$$
where $\mathscr{D}(T_\xi)=H^1_0(\mathcal{D}) \cap H^2(\mathcal{D})$
turns out to be independent of $\xi$. Since $\mathcal{D}$ is
connected, the (elliptic) maximum principle
 implies that $(0-T_\xi)^{-1}$ is
irreducible (see  \cite[Theorem 3.3.5]{davies}   or \cite[Section
11.2]{arendt}). Actually, since $T_\xi$ is the generator of a
holomorphic semigroup, this implies that $(0-T_\xi)^{-1}$ is
positivity improving (see  \cite[p. 306]{nagel}),  i.e.
$(0-T_\xi)^{-1}\varrho(x)
> 0$ for almost every $x \in \mathcal{D}$ provided $\varrho \in
L^2(\mathcal{D})$, $\varrho(x) \geq 0$ for almost every $x \in
\mathcal{D}$ and $\varrho \neq 0$. Now, let $\psi \in X_2$,
$\psi(x,\xi) \geq 0$ for almost every $(x,\xi) \in \mathcal{D}
\times E$, $\psi \neq 0$. Then, $\mathcal{K}_f \psi \geq 0$ and
$(0-T_\xi)^{-1}\mathcal{K}_f\psi(x,\xi) >0$ for almost every
$(x,\xi) \in \mathcal{D} \times E$. It is easy to see that this
exactly means that $(0-\mathcal{T})^{-1}\mathcal{K}_f \psi(x,\xi)
>0$ for almost every $(x,\xi) \in \mathcal{D}\times E$ and the
irreducibility of $(0-\mathcal{T})^{-1}\mathcal{K}_f$ follows. Since
$(0-\mathcal{T})^{-1}\mathcal{K}$ is a compact operator by Theorem
\ref{compacite}, the conclusion follows from Theorems \ref{exist}
and \ref{supinfabs}.
\end{proof}
\subsection{Explicit bounds}

In this section, we derive explicit bounds for the effective
multiplication factor $\ke$. As we did in Section \ref{lowtran}, the
strategy consists in applying Theorem \ref{caracdiff} to suitable
test-functions. We assume the hypothesis of Theorem \ref{prindiff}
to be met. Moreover, we assume here that the diffusion coefficient
$\mathbb{D}(\cdot,\cdot)$ is \textit{degenerate}, i.e.
$$\mathbb{D}(x,\xi)=\mathbb{D}_0(x)d_1(\xi),\qquad \qquad (x,\xi) \in \mathcal{D} \times E,$$
where $\mathbb{D}_0(\cdot)$  is a matrix-valued application
satisfying the ellipticity condition \eqref{ellip},
$\mathbb{D}(\cdot) \in W^{1,2}_\mathrm{loc}(\mathcal{D})$  and
$d_1(\cdot)$ is a\textit{ bounded }real-valued application with
$$\operatornamewithlimits{ess\,inf}_{\xi \in E}d_1(\xi) > 0.$$ Let $\lambda_0$ be the principal eigenvalue of the following
elliptic problem in $L^2(\mathcal{D})$
\begin{equation}\label{diff0}\begin{cases}
\di(\mathbb{D}_0(x)\nabla \varrho(x))+\lambda_0 \varrho(x)=0,\qquad \qquad(x \in \mathcal{D})\\
\varrho_{|\partial \mathcal{D}}(x)=0 \qquad \qquad (x \in \partial
\mathcal{D}).\end{cases}\end{equation}

It is well-known \cite{davies} that $\lambda_0 > 0$  and that there
exists a \textit{positive} eigenfunction $\varrho_0$ solution to
\eqref{diff0}. Set $\mathcal{E}^+=\{\psi \in
L^2(E,\d\xi)\,;\,\psi(\xi)
>0 \:  \text{a.e. } \xi \in E\}.$ In the spirit of Section
\ref{lowtran}, for almost every $x \in \mathcal{D}$, define
$\mathcal{K}^{\lambda_0}_f(x)$ as the following operator on
$L^2(E,\d\xi)$:
$$\mathcal{K}_{f}^{\lambda_0} (x)\::\psi \in L^2(E,\d\xi) \mapsto
\int_{E}\dfrac{\Sigma_f(x,\xi,\xi')}{\lambda_0d_1(\xi)+\sigma(x,\xi)}\psi(\xi')\d\xi'
\in L^2(E,\d\xi)$$ and let  $I_f$ be the set of all $\beta \geq 0$
for which there exists $\psi \in \mathcal{E}^+$ such that
$$[\mathcal{K}_f^{\lambda_0}(x)\psi](\xi) \geq \beta \psi(\xi),\quad
\text{ for almost every } \quad (x,\xi) \in \mathcal{D} \times E.$$
\begin{proposition}\label{lowdiff2} Setting ${\beta}_0:=\sup\{\beta, \beta \in
I\}$, one has $\ke \geq \beta_0$. In particular, a necessary
condition to the reactor to be sub-critical is $\beta_0 < 1$.
\end{proposition}
\begin{proof} As in Section \ref{lowtran}, the set $I$ is closed. Therefore,
there exists ${\psi} \in \mathcal{E}^+$ such that
$[\mathcal{K}_f^{\lambda_0}(x)\psi](\xi) \geq \beta_0 {\psi}(\xi)$
for almost every $(x,\xi) \in \mathcal{D} \times E.$ Now, set
$\varphi_0(x,\xi)=\varrho_0(x)\psi(\xi),$ then, $\varphi \in W_2^+$
and
\begin{equation*}\begin{split}
-\di(\mathbb{D}(x,\xi)\nabla_x \varphi(x,\xi))&=-d_1(\xi)\psi(\xi)\di(\mathbb{D}_0(x)\nabla \varrho_0(x))\\
&=\lambda_0d_1(\xi)\psi(\xi)\varrho_0(x)\qquad \qquad (x,\xi) \in
\mathcal{D} \times E.\end{split}\end{equation*} Consequently, thanks
to \eqref{caracdiff} one has
\begin{equation*}\begin{split}\displaystyle \frac{1}{\ke} &\leq
\underset{(x,\xi) \in \mathcal{D} \times E}
{\operatornamewithlimits{ess\,sup}}\frac{[\lambda_0d_1(\xi)+\sigma(x,\xi)]\varrho_0(x)\psi(\xi)-\varrho_0(x)\displaystyle\int_E\Sigma_s(x,\xi,\xi')\psi(\xi')\d\xi'}{\varrho_0(x)
\displaystyle\int_E\Sigma_f(x,\xi,\xi')\psi(\xi')\d\xi'}\\
&\leq \underset{(x,\xi) \in \mathcal{D} \times E}
{\operatornamewithlimits{ess\,sup}}\frac{[\lambda_0d_1(\xi)+\sigma(x,\xi)]\psi(\xi)}{
\displaystyle\int_E\Sigma_f(x,\xi,\xi')\psi(\xi')\d\xi'}=
\underset{(x,\xi) \in \mathcal{D} \times E}
{\operatornamewithlimits{ess\,sup}}\frac{\psi(\xi)}{[\mathcal{K}_f^{\lambda_0}(x)\psi](\xi)}
\end{split}\end{equation*}
which proves that $\frac{1}{\ke} \leq \frac{1}{\beta_0}.$
\end{proof}In the same spirit, for almost every $x \in \mathcal{D}$, define
$\mathcal{K}^{\lambda_0}(x)$ as the operator on $L^2(E,\d\xi)$ given
by
$$\mathcal{K}^{\lambda_0} (x)\::\psi \in L^2(E,\d\xi) \mapsto
\int_{E}\dfrac{\Sigma(x,\xi,\xi')}{\lambda_0d_1(\xi)+\sigma(x,\xi)}\psi(\xi')\d\xi'
\in L^2(E,\d\xi)$$ where we set
$\Sigma(x,\xi,\xi')=\Sigma_s(x,\xi,\xi')+\Sigma_f(x,\xi,\xi')$. As
in Section \ref{lowtran}, set
$$\overline{\vartheta}:=\inf_{\psi \in \mathcal{E}^+}\esup_{(x,\xi) \in \mathcal{D} \times E}
\dfrac{[\mathcal{K}^{\lambda_0}  (x)\psi](\xi)}{\psi(\xi)}, \quad
\text{ and } \quad \underline{\vartheta}:=\sup_{\psi \in
\mathcal{E}^+}\esinf_{(x,\xi) \in \mathcal{D} \times E}
\dfrac{[\mathcal{K}^{\lambda_0}  (x)\psi](\xi)}{\psi(\xi)}.$$
 Then, one has the following
bounds of $\ge$, in the spirit of Propositions \ref{estimate1} \&
\ref{estimate1bis}.
\begin{proposition}\label{estimate12}   Under the assumptions of Theorem \ref{prindiff}, if $\underline{\vartheta}
> 1$, then $\ge \geq \underline{\vartheta}.$ On the other hand, if $\overline{\vartheta}
> 1$, then $\ge \leq \overline{\vartheta}.$
\end{proposition}
\begin{proof} The proof is very similar to that of Prop.
\ref{estimate1} \& \ref{estimate1bis}. We only prove the first part
of the result, the second part proceeding along the same lines.
Assume thus that $\underline{\vartheta} >1$. For any $\vartheta \in
(1,\underline{\vartheta})$, let $\psi_0 \in \mathcal{E}^+$ be such
that $\esinf_{(x,\xi) \in \mathcal{D} \times E}
\frac{[\mathcal{K}^{\lambda_0} (x)\psi_0](\xi)}{\psi_0(\xi)} \geq
\vartheta.$ Choose then the test-function $\varphi (x,
\xi)=\varrho_0(x)\psi_0(\xi)$. Such an application $\varphi$ belongs
to $W^+_2$ and, as in the above proof, for any $\gamma >0$, one sees
that
\begin{multline*}
(\mathcal{T}+\mathcal{K}(\gamma))\varphi (x,\xi)
 =\frac{\varrho_0(x)}{\vartheta}\left(\lambda_0d_1(\xi)+\sigma(x,\xi)\right)
 \bigg(-\vartheta \psi_0(\xi)+[\mathcal{K}^{\lambda_0} (x)\psi_0](\xi) + \\
 (\vartheta-1)[\mathcal{K}_s^{\lambda_0}(x)\psi_0](\xi)
 +\frac{\vartheta-\gamma }{\gamma}[\mathcal{K}^{\lambda_0}_f( x)\psi_0](\xi)\bigg)\end{multline*}
 where $\mathcal{K}^{\lambda_0}_f(x)$ has been already defined and the definition of $\mathcal{K}^{\lambda_0}_s(x)$ is similar
 ($\Sigma_s$ replacing $\Sigma_f$). Then,
from the positivity of $\Sigma_s$ and $\Sigma_f$, the assumption
$\vartheta
> 1$ implies that
$ (\mathcal{T}+\mathcal{K}(\gamma))\varphi$ is nonnegative for any
$\gamma \leq \vartheta$. Consequently, $\tau_+(\varphi) \geq
\vartheta$ and Prop. \ref{tau} implies that $\ke \geq \vartheta$.
Since $\vartheta \in (1,\underline{\vartheta})$ is arbitrary, one
obtains $\ke \geq \underline{\vartheta}$.\end{proof}

Whenever $E$ is of finite Lebesgue measure, one has the following
practical criteria, already stated by C. V. Pao \cite[Theorem
5.3]{pao2} using completely different arguments.
\begin{corollary}
Assume $E$ to be of finite Lebesgue measure. If
\begin{equation}\label{cvpao}\lambda_0\,d_1(\xi)+\sigma(x,\xi) <
\int_{E}[\Sigma_s(x,\xi,\xi')+\Sigma_f(x,\xi,\xi')]\d\xi'\:\:\:\:\:\:\:\:(x,\xi)
\in \mathcal{D} \times E,
\end{equation}  then, the reactor core is non super-critical, i.e. $\ke \geq 1$.\end{corollary}
\begin{proof} Since $E$ is of finite Lebesgue measure, the constant function
$\psi=\mathbf{1}_E$ such that $\psi(\xi)=1$  for any $\xi \in E$
belongs to $\mathcal{E}^+$. Then, assumption \eqref{cvpao} means
exactly that $[\mathcal{K}^{\lambda_0}(x)\mathbf{1}_E](\xi)
>\mathbf{1}_E(\xi)$ for almost any $(x,\xi) \in \mathcal{D} \times E$.
Therefore, $\underline{\vartheta} >1$ and the conclusion follows
from Prop. \ref{estimate12}.\end{proof}
\begin{remark} Notice that, under the above assumption, one has $$\ge
\geq \esinf_{(x,\xi)}
\frac{\int_{E}[\Sigma_s(x,\xi,\xi')+\Sigma_f(x,\xi,\xi')]\d\xi'}{\lambda_0d_1(\xi)+\sigma(x,\xi)}.$$
\end{remark}
\medskip

\noindent \textbf{Acknowledgment.} The author is grateful to Luisa
Arlotti for a careful reading of a first draft of the manuscript and
her valuable suggestions.
%
%
%
%
%

%
%

\end{document}